\documentclass{article}
\usepackage[utf8]{inputenc}

\usepackage{mathrsfs, amssymb}
\usepackage{amsthm}
\usepackage[intlimits]{empheq}
\usepackage{bbm} 
\usepackage{amsmath}
\usepackage{xcolor}
\usepackage[scientific-notation=true]{siunitx}
\usepackage{float}
\usepackage{dsfont}	
\usepackage[nottoc,numbib]{tocbibind}
\usepackage{bbm}
\usepackage{subcaption}
\usepackage[inner=2cm,outer=2cm,top=2cm,bottom=2cm]{geometry}
\usepackage[ruled,vlined]{algorithm2e}
\usepackage{enumitem}
\usepackage{hyperref}
\usepackage[multiple]{footmisc}
\usepackage[parfill]{parskip}
\usepackage{authblk}
\usepackage{cleveref}
\usepackage[maxcitenames=8,backend=bibtex,style=trad-abbrv,sorting=nyt,sortcites]{biblatex} %style=numeric
\makeatletter
\def\blx@mksc@init{%
  \blx@mkcp@init
  \def\blx@mkcp@nil{\noexpand\blx@mkcp@nil\noexpand}%
  \def\i{\blx@mkcp@nil\i}\def\j{\blx@mkcp@nil\j}%
  \def\do##1{%
    \ifx##1\relax
    \else
      \def##1{\blx@mkcp@nil##1}%
      \expandafter\do
    \fi}%
  \expandafter\do\@uclclist\relax
  \let\(=$\let\)=$}

\def\blx@mksc@eval{%
  \ifx\@let@token\blx@mksc@end
    \expandafter\blx@mksc@end
  \fi
  \ifx\@let@token\bgroup
    \expandafter\blx@mksc@group
  \fi
  \ifx\@let@token\@sptoken
    \expandafter\blx@mksc@space
  \fi
  \ifx\@let@token\blx@mkcp@nil
    \expandafter\blx@mksc@getone
  \fi
  \ifx\@let@token\blx@mkcp@iec
    \expandafter\blx@mksc@getiec
  \fi
  \ifx\@let@token\blx@mkcp@bbl
    \expandafter\blx@mksc@gettwo
  \fi
  \ifx\@let@token\blx@mkcp@sgl
    \expandafter\blx@mksc@gettwo
  \fi
  \ifx\@let@token\blx@mkcp@dbl
    \expandafter\blx@mksc@getthree
  \fi
  \ifx\@let@token$%
    \expandafter\blx@mksc@getmath
  \fi
  \if\noexpand\@let@token\relax
    \expandafter\blx@mksc@cs
  \fi
  \blx@mksc@other}

\def\blx@mksc@getmath#1\blx@mksc@other$#2${\blx@mksc@other{{$#2$}}}
\makeatother

\let\P\relax\DeclareMathOperator{\P}{\mathbb{P}}
\DeclareMathOperator{\Q}{\mathbb{Q}}
\DeclareMathOperator{\E}{\mathbb{E}}
\DeclareMathOperator{\R}{\mathbb{R}}
\DeclareMathOperator{\Z}{\mathcal{Z}}
\DeclareMathOperator{\W}{\mathcal{W}}
\DeclareMathOperator{\U}{\mathcal{U}}
\DeclareMathOperator{\KL}{\operatorname{KL}}
\DeclareMathOperator{\Var}{\operatorname{Var}}

\theoremstyle{plain}
\newtheorem{theorem}{Theorem}[section]
\newtheorem{lemma}[theorem]{Lemma}
\newtheorem{corollary}[theorem]{Corollary}
\newtheorem{proposition}[theorem]{Proposition}
\newtheorem{example}[theorem]{Example}

\theoremstyle{remark}
\newtheorem{remark}[theorem]{Remark}

\title{Nonasymptotic bounds for suboptimal importance sampling}
\author[1]{Carsten Hartmann}
\author[1,2]{Lorenz Richter}
\date{\today}

\affil[1]{Institute of Mathematics, BTU Cottbus-Senftenberg, 03046 Cottbus, Germany, \href{mailto:carsten.hartmann@b-tu.de}{carsten.hartmann@b-tu.de}}
\affil[2]{Institute of Mathematics, Freie Universit\"at Berlin, 14195 Berlin, Germany, \href{mailto:lorenz.richter@fu-berlin.de}{lorenz.richter@fu-berlin.de}}

\begin{document}

\maketitle

\begin{abstract}
Importance sampling is a popular variance reduction method for Monte Carlo estimation, where a notorious question is how to design good proposal distributions. While in most cases optimal (zero-variance) estimators are theoretically possible, in practice only suboptimal proposal distributions are available and it can often be observed numerically that those can reduce statistical performance significantly, leading to large relative errors and therefore counteracting the original intention. In this article, we provide nonasymptotic lower and upper bounds on the relative error in importance sampling that depend on the deviation of the actual proposal from optimality, and we thus identify potential robustness issues that importance sampling may have, especially in high dimensions. We focus on path sampling problems for diffusion processes,  for which generating good proposals comes with additional technical challenges, and we provide numerous numerical examples that support our findings.
\end{abstract}

\section{Introduction}

The numerical approximation of expectations by the Monte Carlo method is ubiquitous in various disciplines such as quantitative finance \cite{glasserman2013monte, glasserman2005importance}, machine learning \cite{bishop2006pattern}, computational statistics \cite{gelman1998simulating} or statistical physics \cite{stoltz2010free}, to name just a few. Depending on the problem at hand, this estimation problem can be more or less difficult, but it turns out that a major challenge are potentially large statistical errors of naive sampling strategies. It is therefore a common goal to build estimators that have a small variance, as compared to the quantity of interest, and thus a small relative error. A typical situation, in which variance reduction is indispensable, is the simulation of rare events with its characteristic exponential divergence of the relative error with the parameter that controls the rarity of the quantity of interest (e.g. a level when computing level-crossing probabilities).

There are multiple strategies for variance reduction in Monte Carlo estimation \cite{asmussen2007stochastic}. In this article, we focus on importance sampling. The idea here is to sample from an alternative probability measure and reweight the resulting random variables with the likelihood ratio in order to produce an unbiased estimator for the quantity of interest. Naturally, the question arises which probability distribution to choose. In theory, under appropriate assumptions, there exists an optimal proposal that yields a zero-variance estimator and therefore removes all the stochasticity from the problem. However this measure depends on the quantity of interest and is therefore practically useless. Coming up with feasible proposals on the other hand is a science in itself, and various numerical experiments demonstrate that it is indeed a crucial one, as making bad choices can even increase the relative error of importance sampling estimators significantly, and therefore counteract the original intention. 
Loosely speaking, importance sampling gets increasingly difficult and sensitive to small deviations from an optimal proposal distribution if the quantity of interest is mainly supported on small regions which have little overlap with the regions of the proposal measure; such a phenomenon is more likely to appear in high dimensions. Moreover, concentration of measure, that may lead to degeneracies of likelihood ratios when the probability of certain events becomes exponentially small, is more likely to occur in high dimensions \cite{polyak2017does, bengtsson2008curse}. 

To better understand the robustness (or better: fragility) of the optimal proposal in importance sampling is the main goal of this article. In applications, one often faces situations that the probability measures admit densities on a subset of  $\R^d$ or a function space like the space of (semi-)continuous trajectories with values in $\R^d$ (called: path space).  In this article, we shall put special emphasis on the latter case, specifically on diffusion processes that are particularly relevant e.g. in molecular dynamics \cite{hartmann2012efficient}, mathematical finance \cite{glasserman2013monte}, or climate modelling \cite{Ragone2018}; what the aforementioned examples have in common, is that the quantities of interest are often related to rare events or large deviations from a mean or an equilibrium state, and, often, the dynamics exhibits metastability, i.e. it features rare transitions between semi-stable equilibria. To simulate these systems, variance reduction techniques like importance sampling are indispensable, and we will provide quantitative bounds on the relative error that explains the fragility of importance sampling in these situations. Some of those bounds are formulated on an abstract measurable space, but they can be readily applied to the density case. For the path space measures, we deduce some additional bounds that, in particular, highlight the challenges due to high dimensionality or long trajectories.

\subsection{Literature overview}

Importance sampling is a classic variance reduction method in Monte Carlo simulation and introductions can be found in many textbooks, such as in \cite[Section 9]{owen2013monte} or \cite{glasserman2013monte, liu2008monte}, however, mostly for the finite-dimensional case $\R^d$. The non-robustness of importance sampling in high dimensions is well known and has often been observed in numerical experiments \cite{bickel2008sharp, snyder2008obstacles, meng1996simulating, glasserman1997counterexamples}. Recently, the authors of \cite{chatterjee2018sample} have proved that the sample size required for importance sampling to be accurate scales exponentially in the KL divergence between the proposal and the target measure, when accuracy is understood in the sense of the $L^1$ error, rather than the commonly used relative error. (Clearly, an unbounded $L^1$ error implies that the relative error will be unbounded.) Similar results can be found in \cite{agapiou2015importance}, in which the authors analyze a self-normalized importance sampling estimator, in connection  with inverse problems and filtering. Necessary conditions that any importance sampling proposal distribution has to satisfy have been derived in  \cite{sanz2018importance}, using the more general $f$-divergences and adopting an information-theoretic perspective.

An important class of techniques for building proposal distributions is known by the name \textit{sequential importance sampling}, where we recommend \cite{doucet2001introduction} for a comprehensive review. Closely related are methods based on interacting particle systems and nonlinear (mean-field) Feynman-Kac semigoups, in which the variance is controlled by adaptively annihilating and generating particles to approximate good proposal distributions \cite{delmoral2013ips}. Adaptive importance sampling for rare events simulation has been pioneered in \cite{dupuis2004importance,dupuis2007subsolutions}; it is typically based on exponential change of measure techniques and the theory of large deviations, dating back to the seminal work \cite{siegmund1976importance}. For diffusion processes, large deviation principles can be used to approximate the optimal change of measure in the small noise regime, where the resulting change of measure turns out to be be asymptotically optimal \cite{vanden2012rare,spiliopoulos2015nonasymptotic}. Pre-asymptotic approximations to the optimal proposal are necessary when studying escape problems, for which the time horizon of the problem is either indefinite or infinitely large, a case that has been analysed in \cite{dupuis2015escaping}. A non-asymptotic variant of the aforementioned approaches for finite noise diffusions is based on the stochastic control formulation of the optimal change of measure \cite{hartmann2017variational,hartmann2019chaos}. Furthermore we should note that there have been many attempts to find good (low-dimensional) proposal by taking advantage of specific structures of the problem at hand, using simplified models that approximate a complicated multiscale system   \cite{dupuis2012importance,spiliopoulos2013large,hartmann2016model,zhang2014optimal}. 
Recently, the scaling properties of certain approximations to control-based importance sampling estimators with the system dimension have been analyzed in \cite{nusken2020solving}, suggesting that the empirical loss function that is used to numerically approximate the optimal proposal distribution is essential.

\subsection{Outline of the paper}

In \Cref{sec: importance sampling bounds divergences} we define importance sampling in an abstract setting and recall the notions of divergences between proposal and target measures, while refining a bound on the relative error and highlighting robustness issues in high dimensions. In \Cref{sec: importance sampling in path space} we move to importance sampling of stochastic processes. We translate the bounds from the previous section to this setting and derive an exact formula for the relative error with which we can state novel bounds that allow for interpretations with respect to robustness in higher dimensions and long time horizons. When focusing on PDE methods in \Cref{sec: PDE methods} we can essentially re-derive bounds from the previous section. In \Cref{sec: small noise diffusions} we comment on how our bounds can help to understand potential issues in the small noise regime. Finally, in \Cref{sec: numerical examples} we present a couple of numerical examples with which we illustrate the previously discussed issues. We conclude the article with \Cref{sec: conclusion} and discuss future perspectives for importance sampling in high dimensions. The article contains an appendix that records some proofs and various technical lemmas. 

\section{Importance sampling bounds based on divergences}
\label{sec: importance sampling bounds divergences}
Let us consider the probability space $(\Omega, \mathcal{F}, \nu)$, on which we want to compute expected values\footnote{As a remark on our notation, let us mention that we sometimes endow the expectation operator with a subscript indicating with respect to which measure the expectation is taken, e.g. $\E_\nu$ indicates that the expectation is considered with respect to the measure $\nu$. When explicitly writing down the corresponding random variable, e.g. $\E[X]$, it is usually clear from the context with respect to which measure the expectation shall be understood, and we omit the subscript.}\footnote{The exponential form, $e^{-\W}$, constrains our observable to be positive. We make this choice in order to be able to have a zero variance proposal density without additional tricks, as the optimal proposal measure $\nu^*$ defined in \eqref{eq: optimal proposal density} has to be non-negative. Assuming strict positivity is convenient in order to get variational dualities that rely on logarithmic transformations, cf. \cite{hartmann2017variational}. An extension of importance sampling to observables with negative parts can for instance be found in \cite{owen2000safe}.}
\begin{equation}
    \label{eq: Z}
    \mathcal{Z} = \E\left[e^{-\W(X)}\right],
\end{equation}
 where $X$ is a random variable taking values in $\Omega$ that is distributed according to the measure $\nu$, and $\W\colon  \Omega \to \R$ is some functional of $X$. Later on we will specify $\Omega$ to be either $\R^d$ or the path space $C([0,T],\R^d)$. 
 
 The idea of importance sampling is to sample instead $\widetilde{X} \in \Omega$ from another distribution $\widetilde{\nu}$ and weight the samples back according to the corresponding likelihood ratio (or Radon-Nikodym derivative), provided that $\nu \ll \widetilde{\nu}$, namely
\begin{equation}
\label{eq: IS for densities}
    \mathcal{Z} =\E\left[e^{-\W(\widetilde{X})}\frac{\mathrm d \nu}{\mathrm d \widetilde{\nu}}(\widetilde{X})\right].
\end{equation}
One notorious intention of importance sampling is the reduction of the variance of the corresponding Monte Carlo estimator
\begin{equation}
\label{eq: Monte Carlo estimator}
    \widehat{\mathcal{Z}}^K = \frac{1}{K} \sum_{k=1}^K e^{-\W(\widetilde{X}^k)}\frac{\mathrm d \nu}{\mathrm d \widetilde{\nu}}(\widetilde{X}^k),
\end{equation}
where $K$ is the sample size and $\widetilde{X}^k$ are i.i.d. samples from $\widetilde{\nu}$. We therefore study the relative error
\begin{equation}
\label{eq: def relative error}
    r(\widetilde{\nu}) = \frac{\sqrt{\Var\left(e^{-\W(\widetilde{X})}\frac{\mathrm d \nu}{\mathrm d \widetilde{\nu}}(\widetilde{X})\right)}}{\mathcal{Z}},
\end{equation}
noting that the true relative error of the estimator \eqref{eq: Monte Carlo estimator} is given by $r(\widetilde{\nu}) / \sqrt{K}$. It can be readily seen that choosing the optimal proposal measure $\widetilde{\nu} = \nu^*$ defined via
\begin{equation}
\label{eq: optimal proposal density}
    \frac{\mathrm d \nu^*}{\mathrm d \nu} = \frac{e^{-\W}}{\mathcal{Z}}
\end{equation}
yields an unbiased zero-variance estimator. Of course, this estimator is usually infeasible in practice, as $\mathcal{Z}$ is just the quantity we are after, and therefore not available. In this article, we study the relative error when using any other absolutely continuous, suboptimal proposal measure $\widetilde{\nu} \neq \nu^*$. It turns out that divergences between those measures are helpful in this analysis and we therefore start by noting the equivalence of the squared relative error and the $\chi^2$ divergence between the actual and the optimal proposal measure.

\begin{lemma}[Equivalence with $\chi^2$ divergence]
\label{lemma: r^2 is chi^2}
Let $\widetilde{\nu}$ be a measure that is absolutely continuous with respect to $\nu$, let $\nu^*$ be the optimal proposal measure as defined in \eqref{eq: optimal proposal density} and let $r(\widetilde{\nu})$ be the relative error as in \eqref{eq: def relative error}. Then
\begin{equation}
    r^2(\widetilde{\nu}) = \chi^2(\nu^* | \widetilde{\nu}).
\end{equation}
\end{lemma}
\begin{proof}
By using the definition of the $\chi^2$ divergence in the first step, we compute
\begin{equation}
    \chi^2(\nu^* | \widetilde{\nu}) = {\E}_{\widetilde{\nu}}\left[\left(\frac{\mathrm d \nu^*}{\mathrm d \widetilde{\nu}}\right)^2 - 1 \right]  = {\E}_{\widetilde{\nu}}\left[\left(\frac{\mathrm d \nu^*}{\mathrm d \widetilde{\nu}}\right)^2\right] - {\E}_{\widetilde{\nu}}\left[\frac{\mathrm d \nu^*}{\mathrm d \widetilde{\nu}} \right]^2 =  {\Var}_{\widetilde{\nu}}\left( \frac{\mathrm d \nu^*}{\mathrm d \widetilde{\nu}}\right) = \frac{1}{\mathcal{Z}^2}\Var\left(e^{-\W(\widetilde{X})}\frac{\mathrm d \nu}{\mathrm d \widetilde{\nu}}(\widetilde{X}) \right) = r^2(\widetilde{\nu}).
\end{equation}
\end{proof}
Motivated by known bounds on the $\chi^2$ divergence, we can formulate our first statement, where we quantify the suboptimality by the Kullback-Leibler divergence between the actual and the optimal proposal measure.
\begin{proposition}[Lower bound on relative error]
\label{prop: r lower bounded with e^KL}
Let $\W : \Omega \to \R$, let $\widetilde{\nu}$ be a measure and let $\nu^*$ be the optimal proposal measure as defined in \eqref{eq: optimal proposal density}, then for the relative error \eqref{eq: IS for densities} it holds
\begin{equation}
\label{eq: relative error lower KL bound}
    r(\widetilde{\nu}) \ge \sqrt{e^{\KL(\nu^* | \widetilde{\nu})} - 1}.
\end{equation}
\end{proposition}
\begin{proof}
With Jensens's inequality we have
\begin{equation}
    \KL(\nu^* | \widetilde{\nu}) = {\E}_{\nu^*}\left[\log \frac{\mathrm d \nu^*}{\mathrm d \widetilde{\nu}} \right] \le \log {\E}_{\nu^*}\left[ \frac{\mathrm d \nu^*}{\mathrm d \widetilde{\nu}} \right].
\end{equation}
Combining this with Lemma \ref{lemma: r^2 is chi^2} yields
\begin{equation}
    r^2(\widetilde{\nu}) =  {\E}_{\widetilde{\nu}}\left[\left(\frac{\mathrm d \nu^*}{\mathrm d \widetilde{\nu}}\right)^2 - 1 \right] =  {\E}_{\nu^*}\left[\frac{\mathrm d \nu^*}{\mathrm d \widetilde{\nu}} - 1 \right] \ge e^{\KL(\nu^* | \widetilde{\nu})} - 1
\end{equation}
and therefore the desired statement.
\end{proof}

\begin{remark}[Bounds on the $\chi^2$ divergence]
In the setting of importance sampling the $\chi^2$ divergence also appears in \cite{chen2005another}. A bound of the $\chi^2$ divergence that is sometimes used is $\chi^2(\nu^* | \widetilde{\nu}) \ge \KL(\nu^* | \widetilde{\nu})$, which is essentially based on $x \le e^{x - 1}$ and therefore yields a less tight bound compared to Proposition \ref{prop: r lower bounded with e^KL}. The exponential bound we use instead can for instance be found in \cite[Theorem 4]{dragomir2000some} and \cite[Proposition 4]{sason2014improved} in a discrete setting; here, a lower bound in terms of the total variation distance is provided as well. \cite{gibbs2002choosing} offers a continuous version and some other helpful relations between divergences. An application of the bound to importance sampling relative errors can be found in \cite{agapiou2015importance} and more analysis with respect to more general $f$-divergences has been done in \cite{sanz2018importance}. The statement should also be compared to the results in \cite{chatterjee2018sample}, where the required sample size of importance sampling is proved to be exponentially large in the KL divergence between the proposal and the target measure.
\end{remark}

\begin{remark}[Cross-entropy method]
Note that the expression $\KL(\nu^*| \widetilde{\nu})$ appearing in \eqref{eq: relative error lower KL bound} is exactly the quantity that is minimized in the so-called cross-entropy method \cite{deboer2005tutorial, zhang2014applications}, which aims at approximating the optimal importance sampling proposal in a family of reference proposals.
\end{remark}

\begin{remark}[Exponential dependence on the dimension]
We recall that the KL divergence usually gets larger with increasing state space dimension as can for instance be seen by Lemma \ref{lem: KL dimension dependence} in the appendix, implying that importance sampling is especially difficult in high dimensional settings. Another way of noting bad scaling behavior in high dimensions is motivated by \cite[Proposition 5.7]{nusken2020solving}. Assume\footnote{The factorization of the optimal proposal measure $\nu^*$ assumes a factorization of the quantity $e^{-g}$.} 
\begin{equation}
    \widetilde{\nu} = \bigotimes_{i=1}^d \widetilde{\nu}_i, \qquad \nu^* = \bigotimes_{i=1}^d \nu^*_i,
\end{equation}
where each $\widetilde{\nu}_i$, and $\nu^*_i$ respectively, shall be identical for $i \in \{1, \dots, d \}$. Then
\begin{align}
    r^2(\widetilde{\nu}) = \Var_{\widetilde{\nu}}\left(\frac{\mathrm d \nu^*}{\mathrm d \widetilde{\nu}} \right) = \E_{\widetilde{\nu}_i}\left[\left(\frac{\mathrm d \nu^*_i}{\mathrm d \widetilde{\nu}_i}\right)^2 \right]^d -  \E_{\widetilde{\nu}_i}\left[\frac{\mathrm d \nu^*_i}{\mathrm d \widetilde{\nu}_i} \right]^{2d} =  \E_{\widetilde{\nu}_i}\left[\left(\frac{\mathrm d \nu^*_i}{\mathrm d \widetilde{\nu}_i}\right)^2 \right]^d -  1 \ge C^d - 1,
\end{align}
where $C := \E_{\widetilde{\nu}_i}\left[\left(\frac{\nu^*_i}{\widetilde{\nu}_i}\right)^2 \right] > 1$ if $\widetilde{\nu} \neq \nu^*$ due to Jensen's inequality. This can be compared to \cite[Section 5.2.1]{sanz2018importance}, and, to be fair, we should note that also naive sampling, i.e. choosing $\widetilde{\nu} = \nu$, usually leads to an exponential dependency of the relative error on the dimension.
\end{remark}

We have so far constructed a lower bound for the relative error. In order to get an upper bound, let us first state the following version of a generalized Jensen inequality, which will turn out to be helpful and is essentially borrowed from \cite[Theorem 2]{mitroi2011estimating}.

\begin{proposition}[Generalized Jensen inequality]
\label{prop: generalized Jensen}
Let $\lambda$ and $\nu$ be measures on $(\Omega, \mathcal{F})$, let
\begin{equation}
        J(f, \nu, \varphi) := \E_{\nu}\left[f(\varphi)  \right] - f \left( \E_{\nu}\left[ \varphi \right]\right)
\end{equation}
be the normalized Jensen functional, where $f : \R \to \R$ is convex and $\varphi: \Omega \to \R$ is continuous, and let $m = \inf_{E \in \mathcal{F}} \frac{\nu(E)}{\lambda(E)}$, $M = \sup_{E \in \mathcal{F}} \frac{\nu(E)}{\lambda(E)}$. Then
\begin{equation}
\label{eq: refined bounds on relative error}
    m J(f, \lambda, \varphi) \le J(f, \nu, \varphi) \le M J(f, \lambda, \varphi).
\end{equation}
\end{proposition}

\begin{proof}
See \Cref{app: proofs - sec 2}.
\end{proof}

We can now derive an upper bound as well as a tighter lower bound for the relative error.

\begin{proposition}[Refined bounds on relative error]
\label{prop: lower and upper bound of relative error}
Let $\widetilde{\nu}$ be a measure that is absolutely continuous with respect to $\nu$ and let $\nu^*$ be the optimal proposal measure as in \eqref{eq: optimal proposal density}. Let $m$ and $M$ be as defined in \Cref{prop: generalized Jensen} (with the measures $\nu$ and $\lambda$ being replaced by $\widetilde{\nu}$ and $\nu^*$ respectively). Then for the relative error \eqref{eq: def relative error} it holds
\begin{equation}
     \sqrt{e^{m \KL(\widetilde{\nu} | \nu^*) + \KL(\nu^* | \widetilde{\nu})} - 1} \le r(\widetilde{\nu}) \le \sqrt{e^{M \KL(\widetilde{\nu} | \nu^*) + \KL(\nu^* | \widetilde{\nu})} - 1}.
\end{equation}

\end{proposition}
\begin{proof}
Inspired by \cite{sason2015tight} (which focuses on a discrete probability space) we choose $\nu = \nu^*, \lambda = \widetilde{\nu}, \varphi = \frac{ \mathrm d \nu^*}{\mathrm d \widetilde{\nu}}$ and $f(x) = -\log(x)$ for the expressions in \eqref{eq: refined bounds on relative error} in order to get
\begin{align}
    J(f, \nu^*, \varphi) &=- \E_{\nu^*}\left[ \log\left(\frac{\mathrm d \nu^*}{\mathrm d \widetilde{\nu}}\right)\right] +\log\left(\E_{\nu^*}\left[ \frac{\mathrm d \nu^*}{\mathrm d \widetilde{\nu}}\right] \right) = -\KL(\nu^* | \widetilde{\nu}) + \log\left(\chi^2(\nu^*|\widetilde{\nu}) + 1\right), \\
    J(f, \widetilde{\nu}, \varphi) &= -\E_{\widetilde{\nu}}\left[ \log\left(\frac{\mathrm d \nu^*}{\mathrm d \widetilde{\nu}}\right)\right]+\log\left(\E_{\widetilde{\nu}} \left[\frac{\mathrm d \nu^*}{\mathrm d \widetilde{\nu}}\right]\right) = \KL(\widetilde{\nu} | \nu^*).
\end{align}
With \Cref{prop: generalized Jensen} we then get
\begin{equation}
    m \KL(\widetilde{\nu} | \nu^*) + \KL(\nu^* | \widetilde{\nu}) \le \log\left(\chi^2(\nu^*|\widetilde{\nu}) + 1\right) \le M \KL(\widetilde{\nu} | \nu^*) + \KL(\nu^* | \widetilde{\nu})
\end{equation}
and with \Cref{lemma: r^2 is chi^2} our statement follows.
\end{proof}

\begin{remark}
One should note that $m$ and $M$ depend on $\widetilde{\nu}$ and $\nu^*$, respectively, and are hard to compute in practice. We have $m \in [0, 1]$ and $M \in [1, \infty]$ and indeed it is possible to get $m = 0$ or $M = \infty$. The former case brings back the ordinary Jensen inequality and the lower bound from \Cref{prop: lower and upper bound of relative error} is then equivalent to the one from \Cref{prop: r lower bounded with e^KL}. The case $M = \infty$ on the other hand yields a trivial upper bound, for which we provide an illustration in \Cref{ex: M = infty}.
\end{remark}

\begin{example}[Upper bound for relative error]
\label{ex: M = infty}
In order to illustrate the case where the upper bound in \Cref{prop: lower and upper bound of relative error} becomes meaningless, consider for instance the measure $\nu$ on $[1, \infty) \subset \R$ admitting the one-dimensional density $p(x) = \alpha \frac{1}{x^{\alpha + 1}}$ defined for $x \ge 1$. 

This density is special since for $\alpha \le 1$ we have $\E[X] = \infty$, however for $\alpha \in (1, 2)$ it holds $\E[X] < \infty$, whereas still $\E\left[X^2\right] = \infty$ and therefore $J(x \mapsto x^2, \nu, \varphi) = \infty$ for $\varphi(x) = x$. Now \Cref{prop: generalized Jensen} implies that the upper bound also has to be infinity. Let us illustrate this for the particular choice of the measure $\lambda$ admitting the density $q(x) = 2\alpha \frac{1}{x^{2\alpha + 1}}$. For this choice we have $J(x \mapsto x^2, \lambda, x \mapsto x) < \infty$ for $\alpha \in (1, 2)$, however we compute
\begin{align}
    M = \sup_{\substack{a, b \in [1, \infty]\\ a \neq b}} \frac{\int_a^b p(x) \mathrm dx}{\int_a^b q(x) \mathrm dx} \ge \sup_{a \in [1, \infty)} \frac{\int_a^\infty p(x) \mathrm dx}{\int_a^\infty q(x) \mathrm dx} = \sup_{a \in [1, \infty)} \frac{ \frac{1}{a^\alpha}}{\frac{1}{a^{2\alpha}}} = \sup_{a \in [1, \infty)} a^\alpha = \infty.
\end{align}
In fact \Cref{prop: generalized Jensen} implies that one cannot not find any $\lambda$ for which both $J(x \mapsto x^2, \lambda, x \mapsto x)$ and $M$ are finite.
\end{example}

To conclude this section, let us illustrate our bounds by looking at a concrete example using Gaussians on $\Omega = \R^d$ (which should be compared to \cite[Section 6]{meng1996simulating}).

\begin{example}[High-dimensional Gaussians]
\label{ex: high-dim Gaussians}
Suppose we want to compute $\E\left[ e^{-\alpha \cdot X}\right]$, with a given vector $\alpha \in \R^d$, where $X\sim \mathcal{N}(\mu, \Sigma) =: p$ is distributed according to a multidimensional Gaussian with mean $\mu \in \R^d$ and covariance matrix $\Sigma \in \R^{d \times d}$. Then the optimal importance sampling density is given by
\begin{equation}
    p^*(x) = \frac{e^{-\alpha \cdot x}}{\Z} p(x) = \mathcal{N}(\mu - \Sigma \alpha, \Sigma).
\end{equation}
If we however sample from a perturbed version
\begin{equation}
\label{eq: perturbed density}
    \widetilde{p}^\varepsilon :=  \mathcal{N}(\mu - \Sigma(\alpha + \varepsilon), \Sigma)
\end{equation}
with a vector $\varepsilon \in \R^d$, we get the relative error
\begin{equation}
\label{eq: r for epsilon shifted Gaussian}
    r(\widetilde{p}^\varepsilon) = \frac{1}{\mathcal{Z}} \sqrt{\Var\left(e^{-\alpha \cdot \widetilde{X}} \frac{p}{\widetilde{p}^\varepsilon}(\widetilde{X}) \right)}= \sqrt{e^{\varepsilon \cdot \Sigma\varepsilon} - 1}.
\end{equation}
In this particular case, the computations can be compared to the relative error of a log-normally distributed random variable, see \Cref{app: relative error log-normal}. Taking, for instance, $\varepsilon = (\widetilde{\varepsilon}, \cdots, \widetilde{\varepsilon})^\top, \Sigma = \operatorname{diag}(\sigma^2, \cdots, \sigma^2)$ yields
\begin{equation}
    r(\widetilde{p}^\varepsilon) = \sqrt{e^{d \sigma^2\widetilde{\varepsilon}^2} - 1},
\end{equation}
where we see an exponential dependence on the variance $\sigma^2$, the squared suboptimality parameter $\widetilde{\varepsilon}$ and the dimension $d$. This implies that, in order to control the relative error in high dimensions, any suboptimal importance sampling estimator needs about $K=\mathcal{O}(e^{d\sigma^2\tilde{\varepsilon}^2})$ independent realisations to reach convergence. This observation is in agreement with the seminal result of Bengtsson and Bickel \cite{bengtsson2008curse} that any importance sampling estimator for Gaussians ceases to be asymptotically efficient when $\log(K)/d\to 0$ as $K,d\to\infty$ (see also \cite[Thm.~3.1]{bengtsson2005curse}).

For this example, we can also apply the bound from Proposition \ref{prop: r lower bounded with e^KL}, by noting that $\KL(p^*|\widetilde{p}^\varepsilon) = \frac{1}{2}\varepsilon\cdot\Sigma\varepsilon$, and get
\begin{equation}
r(\widetilde{p}^\varepsilon) \ge \sqrt{e^{\frac{1}{2}\varepsilon \cdot \Sigma\varepsilon} - 1}.
\end{equation}
A comparison to the exact quantity \eqref{eq: r for epsilon shifted Gaussian} reveals that this lower bound is not tight. For an application of \Cref{prop: lower and upper bound of relative error} we note that also $\KL(\widetilde{p}^\varepsilon|p^*) = \frac{1}{2}\varepsilon\cdot\Sigma\varepsilon$, however $m$ and $M$ are intractable. Still, it is intuitively clear that $m$ becomes smaller and $M$ larger, the more the two Gaussians are apart from each other.

We made the particular choice of $\widetilde{p}^\varepsilon$ in \eqref{eq: perturbed density} in order to have an analogy to the path measure setting, which we will discuss in the next section. In fact, the added term $\Sigma \varepsilon$ in \eqref{eq: perturbed density} can be compared to a constant control $\sigma\sigma^\top \varepsilon$ in a stochastic process as in \eqref{eq: uncontrolled SDE}, which, as will be seen in \eqref{eqn: bound eps^2 d T}, yields a completely analogous expression for the relative error, noting that standard $d$-dimensional Brownian motion is distributed according to $W_T \sim \mathcal{N}(0, \Sigma)$ with $\Sigma = T I_{d \times d}$.
\end{example}

\section{Importance sampling in path space}
\label{sec: importance sampling in path space}
Up to now we have formulated importance sampling for measures on an abstract probability space and provided some illustrations for densities. Let us now elevate those considerations to solutions of stochastic differential equations (SDEs) of the form
\begin{equation}
\label{eq: uncontrolled SDE}
\mathrm d X_s = b(X_s,s) \, \mathrm ds + \sigma(X_s,s) \, \mathrm dW_s, \qquad X_t = x_{\mathrm{init}},
\end{equation}
on the time interval $s \in [t,T]$, $0 \le t < T < \infty$. Here, $b: \mathbb{R}^d \times [t, T] \to \R^d$ denotes the drift coefficient, $\sigma: \R^d \times [t,T]\to \R^{d\times d}$ the diffusion coefficient, $(W_s)_{t \le s \le T}$ standard $d$-dimensional Brownian motion, and $x_{\mathrm{init}} \in \mathbb{R}^d$ is the (deterministic) initial condition. Our goal is to compute expectations of the form\footnote{Whether with $X$ (and $r, \mathcal{Z}$ correspondingly) we refer to random variables in $\R^d$ or solutions to SDEs should usually be clear from the context.}
\begin{equation}
\label{eqn: observable W}
\mathcal{Z} = \E\left[e^{-\mathcal{W}(X)} \right], \qquad \mathcal{W}(X) = \int_t^T f(X_s, s) \mathrm ds + g(X_T),
\end{equation}
where $f : \R^{d} \times [0, T] \to \R, g : \R^{d} \to \R$ are given functions. We will usually fix the initial time to be $t=0$, i.e. consider the SDE \eqref{eq: uncontrolled SDE} on the interval $[0,T]$. For fixed initial condition $x_{\mathrm{init}} \in \mathbb{R}^d$, let us introduce the path space
\begin{equation}
\mathcal{C} = C_{x_{\mathrm{init}}}([0,T],\mathbb{R}^d) = \left\{ X: [0,T]   \rightarrow \mathbb{R}^d  \,\, \vert \,\, X \; \text{continuous}, \; X_0 = x_{\mathrm{init}} \right\},
\end{equation}
equipped with the supremum norm and the corresponding Borel-$\sigma$-algebra, and denote the set of probability measures on $\mathcal{C}$ by $\mathcal{P}(\mathcal{C})$.\par\bigskip
As in the previous section, the idea of importance sampling is to not sample from the original path measure $\P \in \mathcal{P}(\mathcal{C})$ that corresponds to paths of SDE \eqref{eq: uncontrolled SDE}, but from a different measure $\P^u\in\mathcal{P}(\mathcal{C})$ and weight back accordingly. Just as in \eqref{eq: IS for densities} one then gets an unbiased estimator via
\begin{equation}
\label{importance sampling expectation}
    \Z = \E\left[e^{-\W(X^u)} \frac{\mathrm d \P}{ \mathrm d \P^u}(X^u) \right],
\end{equation}
where the Radon-Nikodym derivative is now given by Girsanov's theorem (see Lemma \ref{lem: Girsanov}) and it turns out that the SDE corresponding to $\P^u$ is just a controlled version of the original one,
\begin{equation}
\label{eq: controlled SDE}
\mathrm d X_s^u = \left(b(X_s^u,s) + \sigma(X_s^u,s) u(X_s^u,s)\right) \mathrm ds + \sigma(X^u_s,s) \, \mathrm dW_s, \qquad X_t^u = x_{\mathrm{init}}.
\end{equation}
We think of $u: \mathbb{R}^d \times [t,T] \to \mathbb{R}^d$ as a control term steering the dynamics and note (as already hinted at by the notation) the correspondence between $u$ and $\P^u$. As before, our quantity of interest is the relative error, which now depends on the control $u$:
\begin{equation}
\label{eq: relative error path space}
    r(u) = \frac{\sqrt{\Var\left(e^{-\W(X^u)} \frac{\mathrm d \P}{ \mathrm d \P^u}(X^u) \right) }}{\Z}.
\end{equation}
Given suitable conditions, there exists $u^* \in \U$ that brings \eqref{eq: relative error path space}, the relative error of the importance sampling estimator, to zero \cite{hartmann2017variational}. It turns out that there are multiple equivalent perspectives on the problem of finding such a $u^*$, for instance by solving either a (high-dimensional) Hamilton-Jacobi-Bellman PDE, a forward-backward SDE, a stochastic optimal control problem, or a conditioning of path measures -- the corresponding details regarding those equivalences can for instance be found in \cite{nusken2020solving}. Let us just relate to the last perspective, which claims that $\mathcal{W}$ induces a \emph{reweighted} path measure $\mathbb{Q}$ on $\mathcal{C}$ via
\begin{equation}\label{eq:reweighted measure}
\frac{\mathrm d \Q}{\mathrm d \P} = \frac{e^{-\mathcal{W}}}{\mathcal{Z}},
\end{equation}
assuming $f$ and $g$ are such that $\mathcal{Z}$ is finite (which we shall tacitly assume from now on). It turns out that $\Q = \P^{u^*}$ and we realize that the above formula is the same as in \eqref{eq: optimal proposal density}.\par\bigskip

Let us now bring an example that shall illustrate why variance reduction methods are indispensable in certain SDE settings.

\begin{example}[Rare events of SDEs]
\label{ex: rare events of SDEs}
Monte Carlo estimation gets particularly challenging when considering rare events. As a prominent example, let us consider the one-dimensional Langevin dynamics
\begin{equation}\label{smolu}
\mathrm dX_s = -\nabla \Psi(X_s) \, \mathrm ds + \sqrt{\eta} \, \mathrm dW_s, \quad X_0 = x,
\end{equation}
with double well potential $\Psi(x) = \kappa(x^2-1)^2, \kappa > 0$, and noise coefficient $\eta > 0$, as illustrated in \Cref{fig: double well rare events}. We suppose that the dynamics starts in the left well and choose a function $g$ such that $e^{-g}$ is concentrated in the right well, e.g. $g(x) = \rho (x - 1)^2, \rho > 0$. We are interested in computing $\E[\exp(-g(X_T))|X_0=x]$ for, say, $x=-1$. 

\begin{figure}[H]
\centering
\includegraphics[width=0.50\linewidth]{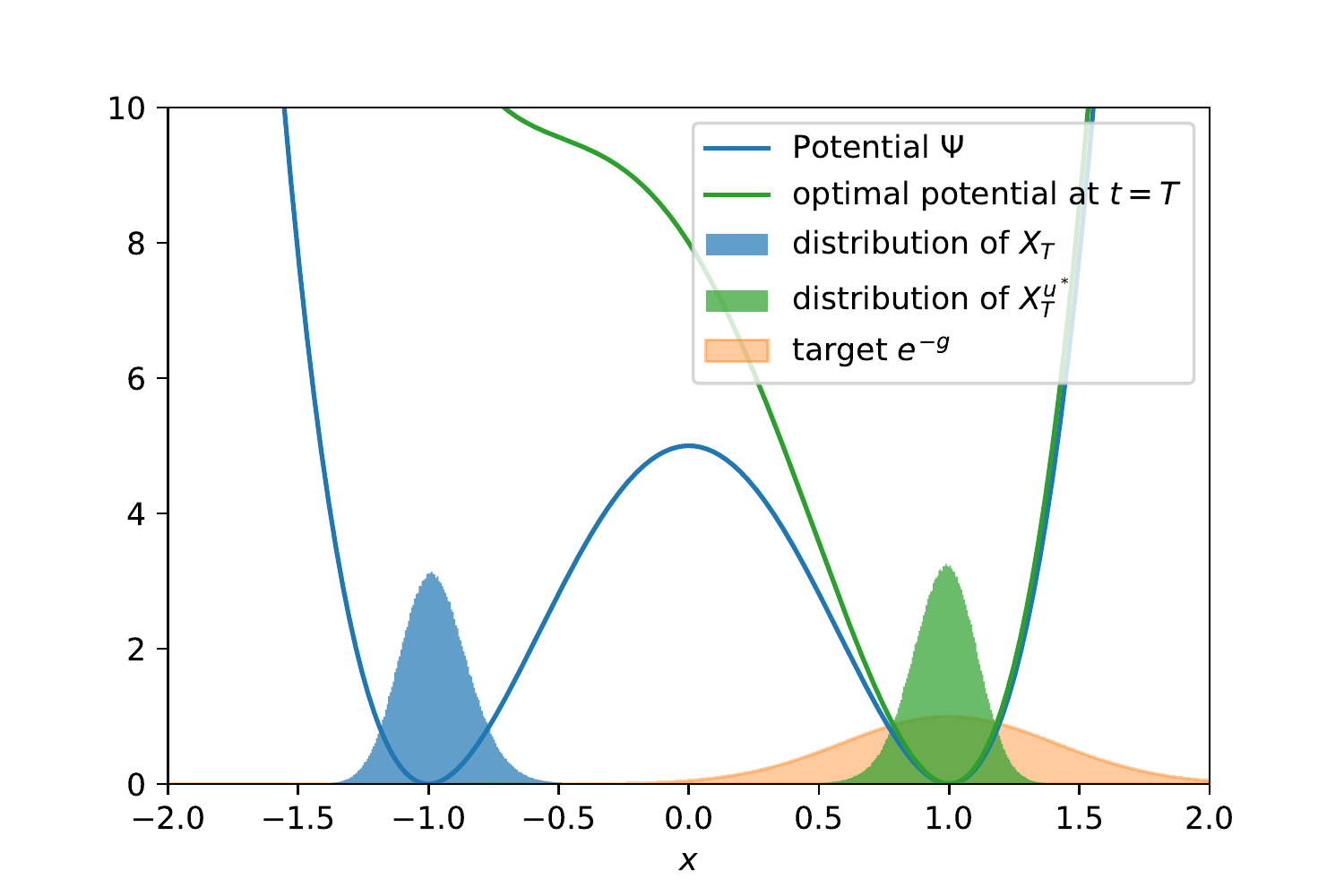}
\caption{Illustration of rare events in a metastable double well potential. We consider the problem described in \Cref{ex: rare events of SDEs} with $\kappa = 5, \rho = 3$ on a time horizon  $T=10$ and display the distributions of $X_T$ as well as $X_T^{u^*}$, which is controlled with the optimal importance sampling control $u^*$ yielding a time-dependent optimal potential.}
\label{fig: double well rare events}
\end{figure}

To understand the difficulties associated with this sampling problem, let $p_T$ be the law of $X_T$ for some $T>0$ and recall that the optimal change of measure is given by the (unnormalized) likelihood $\mathrm d q_T/\mathrm d p_T\propto \exp(-g)$ that is concentrated in the right well. However, regions where $\exp(-g)$ is strongly supported have probability close to zero under $p_T$, for $p_T$ drops to zero quickly for $x>0$. This can be seen as follows: Let $\tau$ be the first exit time of the set $D=\{x\colon x\le 0\}$. By Kramer's law \cite{berglund2011kramers}, the mean first exit time (MFET) satisfies the large deviations asymptotics $\E[\tau]\asymp\exp( 2\Delta\Psi/\eta)$ as $\eta\to 0$, where $\Delta\Psi$ is the energy barrier that the dynamics has to overcome to leave the set $D$, and it turns out that the MFET is independent of the initial condition $x\in D$. Therefore  
\begin{equation}\label{LDP}
    \lim_{\eta\to 0}\eta\log \P(\tau<T) = -2\Delta\Psi\,,\quad T\ll \E[\tau]\,,
\end{equation}
which is is a straight consequence of Kramer's law, combined with the   Donsker-Varadhan large deviations principle that, for a system of the form (\ref{smolu}) states that $\P(\tau<T)\asymp 1-\exp(\lambda_1 T)$ as $T\to\infty$ and $\eta\to 0$, where $\lambda_1\asymp - 1/\E[\tau]$ is the principal eigenvalue of the infinitesimal generator associated with (\ref{smolu}); see, e.g. \cite{Bovier2005}.

Now, by (\ref{LDP}), we can conclude that $p_T(x)\asymp\exp(-2\Delta\Psi/\eta)$ for $x>0$. Since $p_T$ is essentially supported on $(-\infty, 1]$, we can approximate $\exp(-g(x))$ by a step function $\mathbbm{1}_{\{x\in D^c\}}$ on $x\in(-\infty,1]$ and it thus follows that  (up to an exponentially small error)  the relative error for small $\eta>0$ can be approximated by
\begin{subequations}
\begin{align}	
    r(0) & = \sqrt{\frac{\E[\exp(-2g(X_T))]-\E[\exp(-g(X_T))]^2}{\E[\exp(-g(X_T))]^2}}\\
    & \approx \sqrt{\frac{\exp(-2\Delta\Psi/\eta) - \exp(-4\Delta\Psi/\eta)}{\exp(-4\Delta\Psi/\eta)}} 
    \approx \exp(\Delta\Psi/\eta)\,.
\end{align}
\end{subequations}
This kind of exponential behavior is typical for rare event simulation and metastable systems like (\ref{smolu}). So unless our terminal time $T$ is very large or the energy barrier rather small, $X_T$ is usually mostly supported on the left side of the well and therefore does not overlap very much with $e^{-g}$, which leads to an extremely large relative error. Note that this problem gets even more severe with growing values of $\kappa$ and $\rho$.
\end{example}

\subsection{Suboptimal control of stochastic processes and bounds for the relative error}
\label{sec: suboptimal control and RE}
We have stated that, given suitable conditions, there exists $u^* \in \U$ that brings \eqref{eq: relative error path space}, the relative error of the importance sampling estimator, to zero. However, in practice, $u^*$ is usually not available (just as $\nu^*$ is not available in the abstract setting). Let us instead consider the setting where we have the control $u \in \mathcal{U}$ at hand. We want to investigate how the relative error \eqref{eq: relative error path space} behaves depending on how far from optimal $u$ is. For the upcoming analysis, it will turn out that it makes sense to measure the suboptimality and therefore the difference between $\P^u$ and $\P^{u^*}$ in terms of the difference $\delta:=u^* - u$. The first statement is an implication of \Cref{prop: r lower bounded with e^KL}.

\begin{corollary}[Lower bound for relative error on path space]
\label{cor: r bounded by KL paths}
Consider the path measures $\P^u, \P^{u^*} \in \mathcal{P}(\mathcal{C})$ as previously defined and let $\delta=u^* - u$. For the relative error \eqref{eq: relative error path space} it holds
\begin{equation}
    r(u) \ge  \left(\exp\left( \KL(\P^{u^*} | \P^u ) \right) - 1\right)^{\frac{1}{2}} 
    \end{equation}
and therefore
\begin{equation}
\label{eq: lower bound SDE KL}
    r(u) \ge \left(\exp\left( {\E}\left[ \frac{1}{2} \int_0^T |\delta(X_s^{u^*}, s)|^2 \mathrm ds \right] \right) - 1\right)^{\frac{1}{2}}.
\end{equation}
\end{corollary}
\begin{proof}
The first statement is just \Cref{prop: r lower bounded with e^KL} with the abstract measures replaced by path measures. The second statement then follows from Girsanov's theorem as stated in \Cref{lem: Girsanov}.
\end{proof}

One can of course also transfer the more general bound from \Cref{prop: lower and upper bound of relative error} to path measures, however, the computations of the quantities $m$ and $M$ seem even more difficult and impractical than in the density case. In order to still find tighter and more applicable bounds, let us now identify an exact formula for the relative error in the SDE setting.

\begin{proposition}[Formula for path space relative error]
\label{prop: general bound relative error}
Let $X_s^u$ be the solution to SDE \eqref{eq: controlled SDE} and let $\delta=u^* - u$. Then the relative error \eqref{eq: relative error path space} is
\begin{equation}
\label{eqn: general bound relative error - 1}
    r(u) = \left(\E\left[\exp\left(-\int_0^T |\delta(X_s^u, s)|^2 \mathrm ds + 2 \int_0^T \delta(X_s^u, s) \cdot \mathrm dW_s \right) \right] - 1\right)^{\frac{1}{2}},
\end{equation}
or equivalently
\begin{equation}
\label{eqn: general bound relative error - 2}
    r(u) = \left(\E\left[\exp\left(\int_0^T |\delta(X_s^{u + 2\delta}, s)|^2 \mathrm ds \right) \right] - 1\right)^{\frac{1}{2}}.
\end{equation}
\end{proposition}

\begin{proof}
The proof can be found in Appendix \ref{app: proofs}. Alternatively, the second statement follows as well from Proposition \ref{prop: PDE for moment ratio}.
\end{proof}

\begin{remark}
We note that in formula \eqref{eqn: general bound relative error - 1} the forward process is controlled by $u$, whereas in \eqref{eqn: general bound relative error - 2} it is controlled by $u + 2\delta = 2u^* - u$, which of course is usually not available in practice. In the upcoming \Cref{prop: RE bounds independent of x} we will see how we can still make use of the formula.
\end{remark}

\begin{remark}
Note that \Cref{prop: general bound relative error} entails Corollary \ref{cor: r bounded by KL paths} since 
\begin{subequations}
\begin{align}
\label{eqn: change of measure in Jensen bound}
    {\E}_{\P^u}\left[\left(\frac{\mathrm d \P^{u^*}}{\mathrm d \P^u}\right)^2 \right] = {\E}_{\P^{u^*}}\left[\frac{\mathrm d \P^{u^*}}{\mathrm d \P^u} \right] &= \E\left[ \exp\left( \frac{1}{2}\int_0^T |\delta(X_s^{u^*}, s)|^2\mathrm ds +  \int_0^T \delta(X_s^{u^*}, s)\cdot \mathrm d W_s \right)\right] \\
    &\ge \exp\left( \E\left[  \frac{1}{2}\int_0^T |\delta(X_s^{u^*}, s)|^2\mathrm ds \right] \right).
\end{align}
\end{subequations}
Without the change of the measures as in \eqref{eqn: change of measure in Jensen bound} we obtain  
\begin{align}
    {\E}_{\P^u}\left[\left(\frac{\mathrm d \P^{u^*}}{\mathrm d \P^u}\right)^2 \right] \ge \exp\left( \E\left[  -\int_0^T |\delta(X_s^u, s)|^2\mathrm ds \right] \right),
\end{align}
where now the process is controlled by $u$, however this expression has a negative sign in the exponential and is therefore rather useless. The bound
\begin{equation}
    {\E}_{\P^u}\left[\left(\frac{\mathrm d \P^{u^*}}{\mathrm d \P^u}\right)^2 \right] = \E\left[\exp\left(\int_0^T |\delta(X_s^{u + 2\delta}, s)|^2 \mathrm ds \right) \right] \ge \exp\left(\E\left[\int_0^T |\delta(X_s^{u + 2\delta}, s)|^2 \mathrm ds \right]\right) 
\end{equation}
on the other hand, seems more useful. 
\end{remark}

The following corollary derives bounds from the previous \Cref{prop: general bound relative error} that might be useful in practice.

\begin{corollary}[Bounds for path space relative error]
\label{prop: RE bounds independent of x}
Let again $\delta = u^* - u$ and let us assume there exist functions $h_1, h_2: [0, T] \to \R$ such that 
\begin{equation}
    h_1(t) \le |\delta(x, t)| \le h_2(t) 
\end{equation}
for all $x \in \R^d, t \in [0, T]$, then
\begin{equation}
\label{eq: bound relative error time-independent}
     \left(\exp\left( \int_0^T h_1^2(s) \mathrm ds \right)  - 1\right)^{\frac{1}{2}} \le r(u) \le \left(\exp\left( \int_0^T h_2^2(s) \mathrm ds \right)  - 1\right)^{\frac{1}{2}}.
\end{equation}
In particular, if 
\begin{equation}
    \widetilde{\varepsilon}_1 \le |\delta_i(x, t)| \le \widetilde{\varepsilon}_2
\end{equation}
for all components $i \in \{1, \dots, d\}$ and for all $(x, t) \in \R^d \times [0, T]$ with $\widetilde{\varepsilon}_1, \widetilde{\varepsilon}_2 \in \R$, then
\begin{equation}
\label{eqn: bound eps^2 d T}
\left(e^{d\widetilde{\varepsilon}_1^2T} - 1\right)^{\frac{1}{2}}  \le r(u) \le  \left(e^{d\widetilde{\varepsilon}_2^2T} - 1\right)^{\frac{1}{2}}.
\end{equation}
\end{corollary}
\begin{proof}
Both statements follow directly from equation \eqref{eqn: general bound relative error - 2} in Proposition \ref{prop: general bound relative error} by noting that the dependence on the stochastic process and therefore the expectation disappears if we consider bounds on $\delta$ that do not depend on $x$. Two alternative proofs of the corresponding statements can be found in \Cref{app: proofs}.
\end{proof}

\begin{remark}
Note that bounding the suboptimality $\delta$ for all $x$ can be a strong assumption for practical applications, as often, it might vary substantially in $x$. Still, even those conservative bounds often yield lower bounds that render importance sampling a very challenging endeavor. On the contrary, it seems to be hard to make $x$-dependent bounds on $\delta$ useful due to potentially very complex stochastic dynamics. Let us further note that the bounds in \eqref{eq: bound relative error time-independent} imply that errors made over different points in time accumulate, i.e. it does not matter if they have been made at the beginning or the end of a trajectory and neither can they be compensated at later stages.
\end{remark}

Another upper bound on the relative error can be derived by means of the H\"older inequality.

\begin{proposition}[Another bound for path space relative error]
\label{prop: Hoelder bound}
Let $\delta = u^* - u$. For the relative error \eqref{eq: controlled SDE} it holds
\begin{equation}
     r(u) \le \left({\E}\left[\exp \left((1 + \sqrt{2})^2 \int_0^T |\delta(X_s^u,s)|^2 \mathrm ds  \right) \right]^{\frac{1}{1+\sqrt{2}}} - 1\right)^\frac{1}{2} 
     \end{equation}
\end{proposition}

\begin{proof}
See Appendix \ref{proof: Hoelder bound}.
\end{proof}

\begin{remark}
Some intuition of the quality of this bound can be gained when for instance assuming that $\delta(x, t) = \varepsilon$ with a constant vector $\varepsilon = (\widetilde{\varepsilon}, \dots, \widetilde{\varepsilon})^\top \in \R^d$. Then this bound yields $r(u) \le \left(\exp\left((1+ \sqrt{2})d \widetilde{\varepsilon}^2 T \right) - 1\right)^{\frac{1}{2}}$, which is less tight than the bound \eqref{eqn: bound eps^2 d T} in \Cref{prop: RE bounds independent of x}. Nevertheless the bound is useful in that it only depends on the stochastic process controlled by $u$, which is a known quantity. 
\end{remark}

\subsection{PDE methods for the study of relative errors}
\label{sec: PDE methods}
Another means of studying the relativ error $r(u)$ are partial differential equations (PDEs). We will formulate a PDE for the relative error \eqref{eq: def relative error}, which might be helpful for future analysis and by which we can rederive bounds from the previous section.\par\bigskip 

By a slight generalization of \cite{spiliopoulos2015nonasymptotic}, one can identify a PDE for the $u$-dependent second moment (conditioned on $X^u_t = x$),
\begin{equation}
\label{eq: second moment for SDE}
M_u(x, t) = \E\left[e^{-2\mathcal{W}(X^u)} \left(\frac{\mathrm d \P}{\mathrm d \P^u}(X^u)\right)^2\Bigg|X^u_t = x  \right],
\end{equation}
namely
\begin{subequations}
\begin{align}
\label{eqn: PDE second moment}
    (\partial_t + L - \sigma u(x, t) \cdot \nabla - 2 f(x, t) + |u(x, t)|^2) M_u(x, t) &= 0, \qquad &(x, t) \in  {\R}^d \times [0, T) , \\
    M_u(x, T) &= e^{-2g(x)}, \qquad &x \in {\R}^d,
\end{align}
\end{subequations}
where $L = \frac{1}{2}(\sigma \sigma^\top)(x, t):\nabla^2 + b(x, t) \cdot \nabla$ is the infinitesimal generator associated to the SDE \eqref{eq: uncontrolled SDE}.\par\bigskip

Defining $\delta = u^* - u$, this then immediately leads to the PDE
\begin{subequations}
\begin{align}
\label{eqn: PDE for second moment}
    (\partial_t + L + \sigma (\sigma^\top \nabla V(x, t)+ \delta(x, t)) \cdot \nabla - 2 f(x, t) + |\sigma^\top \nabla V(x, t)+ \delta(x, t)|^2) M_u(x, t) &= 0, \quad &(x, t) \in  {\R}^d \times [0, T), \\
    M_u(x, T) &= e^{-2g(x)}, \quad &x \in {\R}^d,
\end{align}
\end{subequations}
which describes the second moment of suboptimal importance sampling. It can be shown that for $\delta = 0$, i.e. under the optimal control $u = u^*$, we recover indeed the zero-variance property of the corresponding importance sampling estimator, see Proposition \ref{prop: u^* implies zero-variance} in the appendix. In the following statement we construct the PDE that is relevant for the relative error $r(u)$ and re-derive a formula that we have already seen before.

\begin{proposition}[PDE for the relative error]
\label{prop: PDE for moment ratio}
Let $\delta = u^* - u$. We consider the second moment as in \eqref{eq: second moment for SDE} and the conditional expectation $\psi(x, t) = \E\left[e^{-\mathcal{W}(X^u)} \Big|X_t = x  \right]$, then the function $h_u:{\R}^d \times [0, T] \to \R$ defined by
\begin{equation}
     h_u(x, t) = \frac{M_u(x, t)}{\psi^2(x, t)},
\end{equation}
solves the PDE
\begin{subequations}
\label{eq: PDE for relative error}
\begin{align}
    \left(\partial_t + L^{u+2\delta} + |\delta(x, t)|^2\right)h_u(x, t) &= 0, \qquad &(x, t) \in  {\R}^d \times [0, T),\\
    h_u(x, T) &= 1,\qquad &x \in {\R}^d,
\end{align}
\end{subequations}
with $L^{u+2\delta} := L + \sigma (u + 2\delta)\cdot \nabla$. This then implies
\begin{equation}
\label{eqn: moment term h(x,t)}
    h_u(x, t) = \E\left[\exp\left(\int_t^T |\delta(X_s^{u + 2\delta}, s)|^2 \mathrm ds \right) \Bigg| X_t^{u + 2\delta} = x \right].
\end{equation}
\end{proposition}

\begin{proof}
We plug the ansatz
\begin{equation}
    M_u(x, t) = h_u(x, t) \psi^2(x, t) = h_u(x, t) e^{-2V(x, t)}
\end{equation}
into the PDE \eqref{eqn: PDE for second moment}. Noting that
\begin{subequations}
\begin{align}
    (\sigma \sigma^\top):\nabla^2(h_u e^{-2V}) &= (\sigma \sigma^\top):\left(\nabla \left( \nabla h_u e^{-2V} - 2h_u \nabla V e^{-2V} \right)\right) \\
    &= e^{-2V}\left((\sigma \sigma^\top):\nabla^2h_u - 4\,\sigma \sigma^\top\nabla V \cdot \nabla h_u + 4\,h_u|\sigma^\top \nabla V|^2 - 2h_u(\sigma \sigma^\top):\nabla^2 V \right),
\end{align}
\end{subequations}
we get the PDE
\begin{align}
    -2h_u \underbrace{\left(\partial_t V  + LV - \frac{1}{2}|\sigma^\top \nabla V|^2 + f \right)}_{=0} + \partial_t h_u + L h_u - \sigma \sigma^\top \nabla V\cdot \nabla h_u + \sigma \delta \cdot \nabla h_u + |\delta|^2 h_u = 0,
\end{align}
from which the statement follows from the identity $u^* = -\sigma^\top \nabla V$ and a specific Hamilton-Jacobi-Bellman equation that is for instance stated in \cite[Problem 2.2]{nusken2020solving}. The probabilistic representation \eqref{eqn: moment term h(x,t)} follows immediately from the Feynman-Kac formula \cite[Theorem 1.3.17]{pham2009continuous}.
\end{proof}

\begin{remark}
First note that $h_u$ from \Cref{prop: PDE for moment ratio} is related to the relative error \eqref{eq: relative error path space} via $r(u) = \sqrt{h_u(x, 0) - 1}$. On the first glance it looks like the PDE \eqref{eq: PDE for relative error} does not depend on $f$ and $g$. This is of course not true and we should note that the PDE depends on $u^*$, which again depends on $f$ and $g$. Finally, note that with \eqref{eqn: moment term h(x,t)} we recover the result \eqref{eqn: general bound relative error - 2} from Proposition \ref{prop: general bound relative error}.
\end{remark}

\subsection{Small noise diffusions}
\label{sec: small noise diffusions}
A prominent application of importance sampling in stochastic processes can be found in the context of small noise diffusions and rare event simulations (relating to \Cref{ex: rare events of SDEs}, see also \cite{spiliopoulos2013large, spiliopoulos2015nonasymptotic, vanden2012rare, dupuis2012importance}). We model small noises with the smallness parameter $\eta > 0$ by considering the SDEs\footnote{To be consistent with the notation from before, we could hide the smallness parameter $\eta$ in the diffusion coefficient, i.e. $\sigma = \sqrt{\eta}\widetilde{\sigma}$. Then the HJB equation that provides the zero variance control is $(\partial_t + \frac{\eta}{2}(\widetilde{\sigma}\widetilde{\sigma}^\top) :\nabla^2 + b \cdot \nabla) V - \frac{1}{2}|\sqrt{\eta}\widetilde{\sigma}\nabla V|^2 + \frac{1}{\eta} f(x, t) = 0, V(x, T) = \frac{1}{\eta}g(x)$ and the relation $V^\eta = \eta V$ yields HJB equation \eqref{eq: HJB epsilon}.}
\begin{equation}
\mathrm d X_s^\eta = b(X_s^\eta,s) \, \mathrm ds + \sqrt{\eta}\,\widetilde{\sigma}(X_s^\eta,s) \, \mathrm dW_s, \qquad X^\eta_t = x_{\mathrm{init}},
\end{equation}
and we want to compute quantities like
\begin{equation}
    \psi^\eta(x, t) = \E\left[e^{-\frac{1}{\eta} \mathcal{W}(X^\eta)} \Big| X_t^\eta = x \right].
\end{equation}
If $\eta$ gets smaller it becomes harder to estimate $\psi^\eta(x, t)$ via Monte Carlo methods as the variance grows exponentially in $\eta$. To be more precise, by Varadhan's lemma \cite[Theorem 4.3.1]{dembo2009large}, using the quantities 
\begin{equation}
    \gamma_1 := -\lim_{\eta \to 0} \eta \log \E\left[e^{-\frac{1}{\eta}\mathcal{W}(X^\eta)}\right] \qquad \text{and} \qquad \gamma_2 := -\lim_{\eta \to 0} \eta  \log \E\left[e^{-\frac{2}{\eta}\mathcal{W}(X^\eta)}\right],
\end{equation}
one gets for the relative error of the uncontrolled process 
\begin{equation}
    r(0) = \sqrt{e^{\frac{2\gamma_1 - \gamma_2 + o(1)}{\eta}}-1},
\end{equation}
asymptotically as $\eta\to 0$. By Jensen's inequality we have $2\gamma_1>\gamma_2$ unless $\mathcal{W}$ is a.s. constant, but we note that even for $2\gamma_1 = \gamma_2$ the relative error explodes in the limit $\eta \to 0$. Let us again consider a controlled process
\begin{equation}
\mathrm d X_s^{u, \eta} = \left(b(X_s^{u, \eta},s) + \widetilde{\sigma}(X_s^{u, \eta},s)u(X_s^{u, \eta},s)  \right) \, \mathrm ds + \sqrt{\eta}\widetilde{\sigma}(X_s^{u, \eta},s) \, \mathrm dW_s, \qquad X^{u, \eta}_t = x_{\mathrm{init}},
\end{equation}
and realize that the optimal importance sampling control that yields zero variance,
\begin{equation}
    u^* = -\widetilde{\sigma}^\top \nabla V^{\eta}  = \eta \widetilde{\sigma}^\top \nabla \log \psi^\eta,
\end{equation}
can be computed via the HJB equation
\begin{equation}
\label{eq: HJB epsilon}
    \left(\partial_t + \frac{\eta}{2}(\widetilde{\sigma} \widetilde{\sigma}^\top)(x, t) : \nabla^2 +  b(x, t) \cdot \nabla\right) V^\eta(x, t) - \frac{1}{2}|(\widetilde{\sigma}^\top \nabla V^\eta)(x, t)|^2 + f(x, t) = 0, \qquad V^\eta(x, T) = g(x).
\end{equation}
Since solving this PDE is notoriously difficult (especially in high dimensions), various  approximations have been suggested that lead to estimators that enjoy log-efficiency or a vanishing relative error in the regime of a vanishing $\eta$. However, since log-efficient estimators still often perform badly in practice (as for instances discussed in \cite{asmussen2011importance, glasserman1997counterexamples}), in \cite{vanden2012rare} it is suggested to replace $u^*$ by the vanishing viscosity approximation $u^0$ based on the  corresponding HJB equation with $\eta = 0$: 
\begin{equation}
    u^0 = -\widetilde{\sigma}^\top\nabla V^0,
\end{equation}
where $V^0$ is the solution to 
\begin{equation}
\label{eq: HJB deterministic problem}
    (\partial_t + b(x, t) \cdot \nabla) V^0(x, t) - \frac{1}{2}|(\widetilde{\sigma}^\top\nabla V^0)(x, t)|^2 + f(x, t) = 0, \qquad V^0(x, T) = g(x).
\end{equation}
While it can be shown that, given some regularity assumptions on $f$ and $g$, it holds \cite{vanden2012rare}
\begin{equation}
    \lim_{\eta \to 0} r(u^0) = 0,
\end{equation}
a large relative error for a small, but fixed $\eta > 0$ is still possible. In our notation from before, this situation corresponds to choosing $\delta = u^* - u^0$ and Propositions \ref{prop: general bound relative error} and \ref{prop: PDE for moment ratio} show that
\begin{equation}
    r(u^0) = \sqrt{\E\left[\exp\left(\int_0^T |u^* - u^0|^2(X_s^{2u^* - u^0}, s) \mathrm ds \right) \right] - 1}.
\end{equation}
Even though this expression converges to zero as $\eta \to 0 $ provided that $V \to V^0$ and $u^* \to u^0$ \cite{fleming1971stochastic}, we expect an exponential dependence on the time $T$ and the dimension $d$ for any fixed $\eta > 0$ (cf. our numerical experiment in \Cref{sec: numerics small noise diffusion}).

In \cite{fleming1971stochastic} it is proved that 
\begin{equation}
     \nabla V = \nabla  V^0 + \eta \nabla  v_1 + o(\eta),
\end{equation} 
uniformly on all compact subsets of $\R^d\times (0,T)$, where $v_1$ solves the PDE stated in \Cref{app: asymptotic expansion small noise}. As a consequence, we can write
\begin{equation}
    | \nabla V - \nabla  V^0| = |\eta \nabla v_1 + o(\eta)| = \eta| \nabla v_1 + o(1)|
\end{equation}
and
\begin{equation}
   r(u^0) = \E\left[\exp\left(\eta^2 \int_0^T |(\sigma^\top\nabla v_1)(X_s^{2u^* - u_0^*}, s)|^2 \mathrm ds + o(\eta^2)\right) \right].
\end{equation}
Specifically, if there exist constants $C_1, C_2 > 0$ such that  $C_1 < |\nabla v_1(x, t)| < C_2$ for all $(x, t) \in \R^d \times (0, T)$, then the relative error grows exponentially as  
\begin{equation}
\label{eqn: relative error small noise diffusion}
         \sqrt{e^{\eta^2C_1^2 T + o(\eta^2)}-1} \le r(u^*) \le \sqrt{e^{\eta^2C_2^2 T + o(\eta^2)}-1}
\end{equation}
due to \Cref{prop: RE bounds independent of x}. We emphasize, however, that it is not clear under which assumptions this uniform bound can be achieved, given that, in practice, $v_1$ can be strongly $x$-dependent as is illustrated with a numerical example in Section \ref{sec: numerics small noise diffusion}.

\begin{remark}
The above considerations show that the relative error is potentially only small if $\eta$ is (much) smaller than $C_1 \sqrt{T}$. This can be compared to equation (5.3) in \cite{spiliopoulos2015nonasymptotic} and in particular to \cite{dupuis2015escaping}, where a concrete example is constructed for which the second moment can be lower bounded by $e^{-\frac{1}{\eta} C_1 + (T-K)C_2}$ for $C_1, C_2, K > 0$, i.e. the time $T$ and the smallness parameter $\eta$ compete. We illustrate the degeneracy with growing $T$ for a toy example in \Cref{fig: small noise rel error T}.
\end{remark}

\section{Numerical examples}
\label{sec: numerical examples}

In this section we provide numerical examples that shall illustrate some of the formulas and bounds derived in the previous sections. We particularly demonstrate that importance sampling can be very sensitive to small perturbations of the optimal proposal measure. Here we focus on path space measures and provide several examples of importance sampling of diffusions. The code can be found at \url{https://github.com/lorenzrichter/suboptimal-importance-sampling}.

\subsection{Ornstein-Uhlenbeck process}

\begin{figure}[H]
\centering
\includegraphics[width=1.00\linewidth]{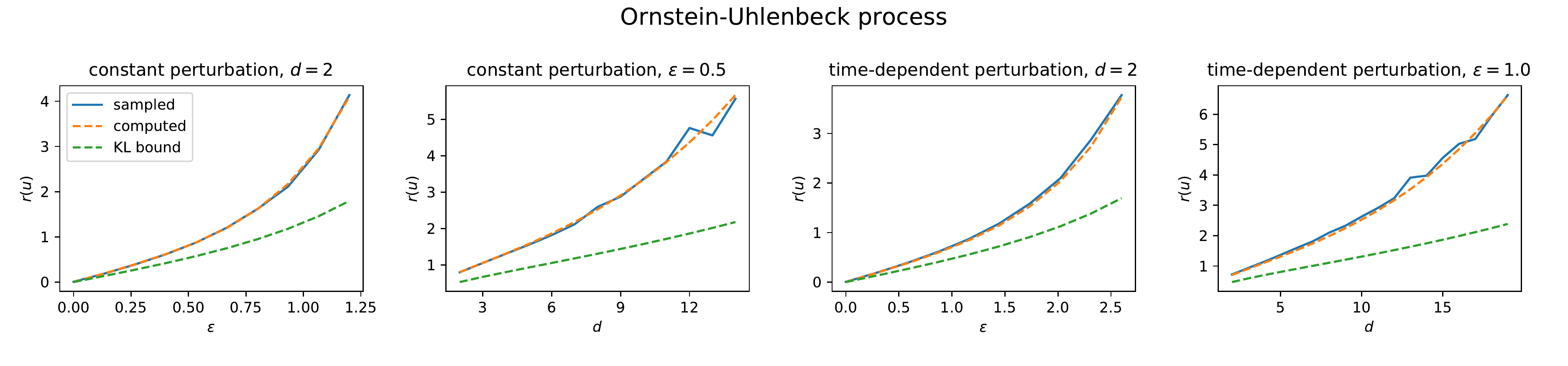}
\caption{Sampled relative error with varying constant or time-dependent perturbation $\varepsilon$ and dimension $d$ compared to the formulas derived in \Cref{prop: RE bounds independent of x} and to the lower bound from \Cref{cor: r bounded by KL paths}.}
\label{fig: OU suboptimal eps and dimensions}
\end{figure}

An example where the optimal importance sampling control is analytically computable is the following. Consider the $d$-dimensional Ornstein-Uhlenbeck process
\begin{equation}
\mathrm dX_s = AX_s  \mathrm d s + B \,\mathrm d W_s, \quad X_0 = 0,
\end{equation}
and its controlled version
\begin{equation}
\mathrm dX_s^u = \left(AX_s^u + B u(X_s^u, s)\right) \mathrm d s + B \,\mathrm d W_s, \quad X_0^u = 0,
\end{equation}
where $A,B \in \R^{d \times d}$ are given matrices. In \eqref{eqn: observable W} we set $f = 0$ and  $g(x) = \alpha \cdot x$, for a fixed vector $\alpha \in \R^d$, i.e. we want to estimate the quantity
\begin{equation}
    \mathcal{Z} = \E\left[e^{-\alpha \cdot X_T}\right].
\end{equation}
As shown in \cite{nusken2020solving}, the zero-variance importance sampling control is given by
\begin{equation}
u^*(x, t) = -B^\top e^{A^\top (T-t)}\alpha.
\end{equation}
We choose $A = -3\, I_{d \times d} + (\xi_{ij})_{1\le i,j \le d}$ and  $B = I_{d \times d} + (\xi_{ij})_{1\le i,j \le d}$, where $\xi_{ij} \sim \mathcal{N}(0, \sigma^2)$ are i.i.d. random coefficients that are held fixed throughout the simulation. We set $T = 1, \sigma = 1$, $\alpha = (1, \dots, 1)^\top$ and first consider the perturbed control
\begin{equation}
    u = u^* + (\varepsilon, \dots, \varepsilon)^\top.
\end{equation}
In the two left panels of \Cref{fig: OU suboptimal eps and dimensions} we display a Monte Carlo estimation of the relative error \eqref{eq: relative error path space} using $K = 10^6$ samples and compare it to the formulas from \Cref{prop: RE bounds independent of x} and the bound from \Cref{cor: r bounded by KL paths}, once with varying perturbation strength $\varepsilon$, once with varying dimension $d$. We see that in both cases the simulations agree with our formula, even though for moderate to large deviations from optimality the estimated values of $r$ are observed to fluctuate.

Let us now look at an example with a time-dependent perturbation of the optimal control. More specifically,  we consider a perturbation that is active only for a certain amount of time $s < T$, namely
\begin{equation}
    u(x, t) = u^*(x, t) + (\varepsilon, \dots, \varepsilon)^\top \mathbbm{1}_{[0, s]}(t),
\end{equation}
where in our experiment we choose $s = 0.2$. In the two right panels of \Cref{fig: OU suboptimal eps and dimensions} we display the same comparisons as before, however now using formula \eqref{eq: bound relative error time-independent} in order to account for the time-dependent nature of the perturbation.

\subsection{Double well potential}
For strongly metastable systems, Monte Carlo estimation is notoriously difficult and variance reduction methods are often indispensable. Importance sampling seems like a method of choice, but we want to illustrate that one has to be very careful with the design of the importance sampling control.

\begin{figure}[H]
\centering
\includegraphics[width=0.7\linewidth]{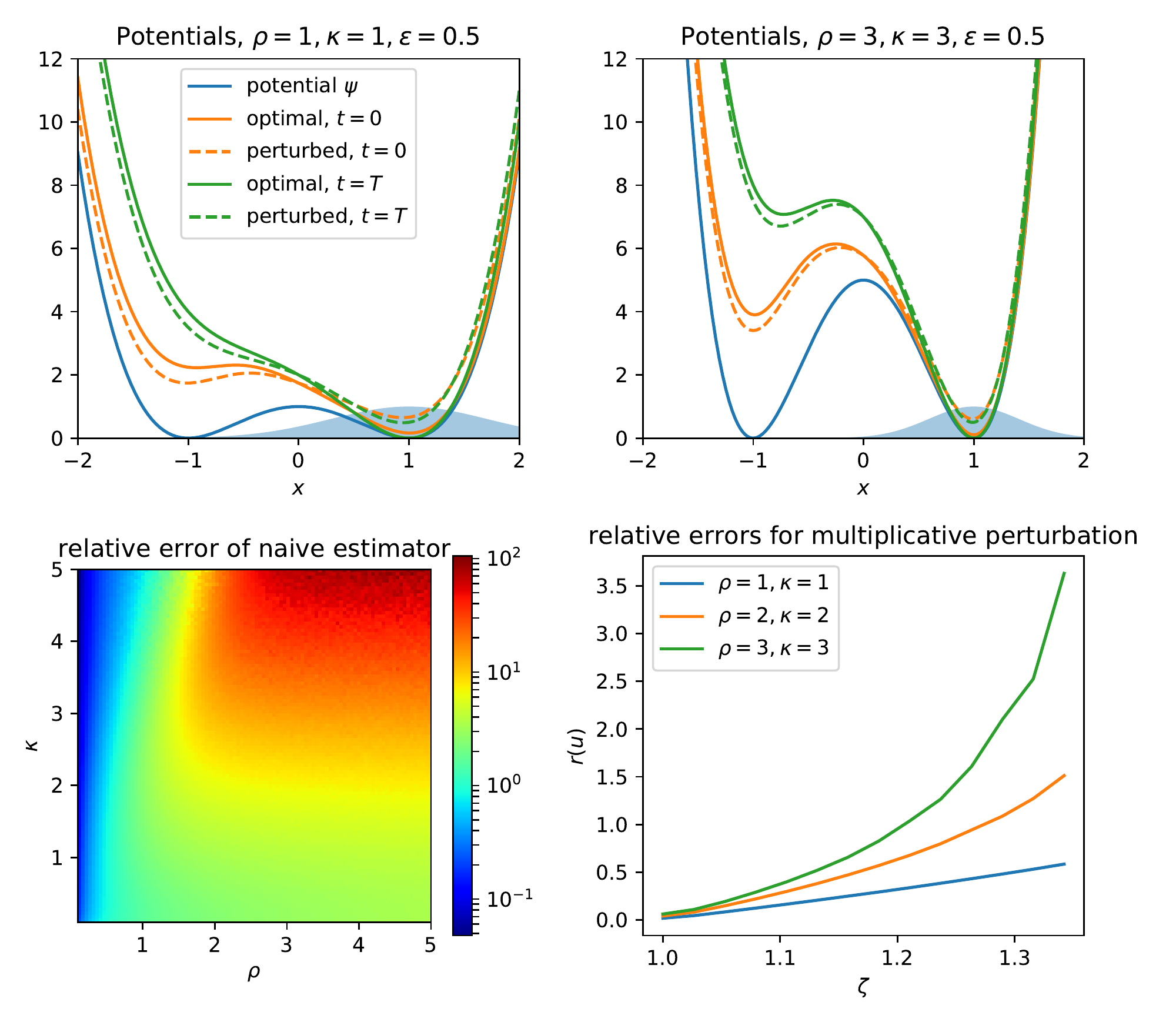}
\caption{Top panels: Double well potentials and optimal tiltings as well as additive perturbations for different values of $\rho$ and $\kappa$. Bottom left: Relative error of the naive Monte Carlo estimator for different values of $\rho$ and $\kappa$. Bottom right: Relative error depending on the multiplicative perturbation factor $\zeta.$}
\label{fig: double well tilted potentials}
\end{figure}

As in \Cref{ex: rare events of SDEs}, let us consider the Langevin SDE
\begin{equation}
\label{SDE_double_well}
\mathrm dX_s = -\nabla \Psi(X_s) \, \mathrm ds + B \, \mathrm dW_s, \quad X_0 = x,
\end{equation}
in $d=1$, where $B \in \R$ is the diffusion coefficient,  $\Psi(x) = \kappa(x^2-1)^2$ is a double well potential with $\kappa > 0$ and $x = -1$ is the initial condition. For the observable in  \eqref{eqn: observable W} we consider $f = 0$ and $g(x) = \rho (x-1)^2$, where $\rho > 0$; the terminal time is set to $T=1$. Note that choosing higher values for $\rho$ and $\kappa$ accentuates the metastable features, making sample-based estimation of $ \E\left[\exp(-g(X_T))\right]$ more challenging. For an illustration, the two top panels of Figure \ref{fig: double well tilted potentials} show the potential $\Psi$ and the weight from \eqref{eq:reweighted measure}, $e^{-g(x)}$, for different values of $\rho$ and $\kappa$ and for $B=1$. We also plot the `optimally tilted potentials' $\Psi^* = \Psi + BB^\top V$, noting that $-\nabla \Psi^* = -\nabla \Psi + Bu^*$. In the bottom left panel we show the relative error of the naive estimator depending on different values of $\rho$ and $\kappa$.

As before, let us perturb the optimal control, this time both in an additive and multiplicative way, namely
\begin{equation}
    u = u^* + \varepsilon = -B^\top \nabla (V - B^{-\top}\varepsilon \cdot x) \qquad\text{and}\qquad u = \zeta u^*,
\end{equation}
where $\varepsilon \in \R^d, \zeta \in \R$ specify the perturbation strengths. In the bottom right panel of \Cref{fig: double well tilted potentials} we show the relative error for the multiplicative perturbation and see that for higher values of $\rho$ and $\kappa$ the exponential divergence becomes more severe, demonstrating that the robustness issues of importance sampling are particularly present in metastable settings.

\begin{figure}[H]
\centering
\includegraphics[width=1.0\linewidth]{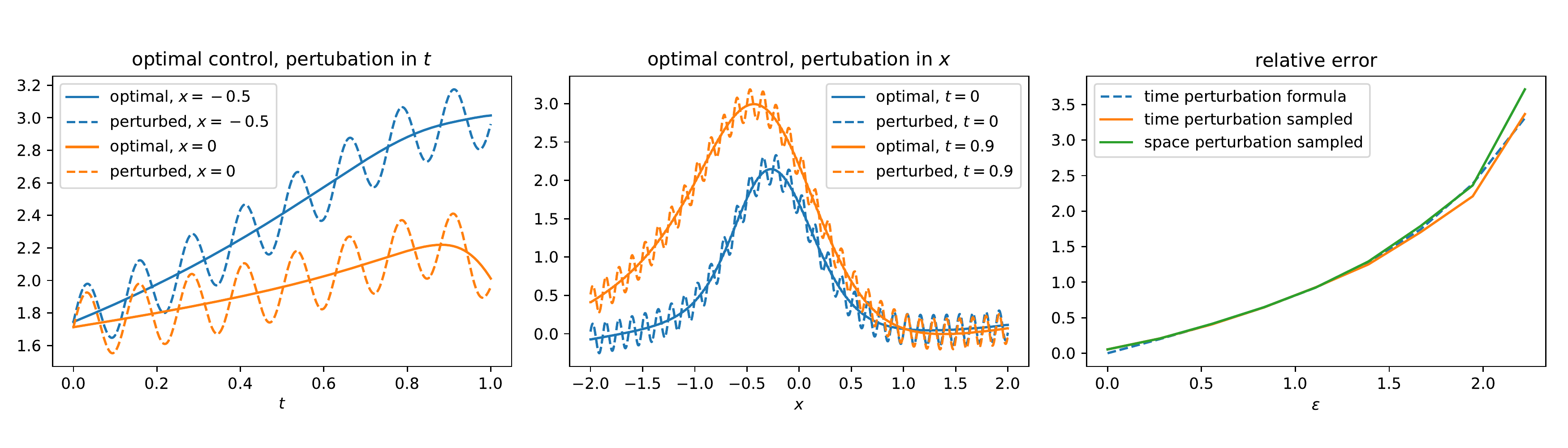}
\caption{Left: Optimal importance sampling control and time perturbation for two different values of $x$. Middle: Optimal importance sampling control and space perturbation for two different values of $t$. Right: Relative error of suboptimal importance sampling estimators depending on the perturbation strength $\varepsilon$; here, the dashed line refers to the exact formula (\ref{eq: relative error double well time-dependent perturbation}).}
\label{fig: double well sine perturbations}
\end{figure}

Let us now consider perturbations depending either on time or space,
\begin{equation}
    u_1(x, t) = u^*(x, t) + \varepsilon \sin(\alpha t)  \qquad \text{and} \qquad u_2(x, t) = u^*(x, t) + \varepsilon \sin(\alpha x),
\end{equation}
as illustrated in Figure \ref{fig: double well sine perturbations} with $\alpha = 50$.

In the former case we can analytically compute the relative error due to \Cref{prop: RE bounds independent of x} to be
\begin{equation}
\label{eq: relative error double well time-dependent perturbation}
    r_1(\varepsilon) = \sqrt{\exp\left(\varepsilon^2\left(\frac{T}{2} - \frac{\sin(2\alpha T)}{4 \alpha} \right)\right) - 1}.
\end{equation}
Let us again illustrate how the relative error depends on the perturbation strength $\varepsilon$. In the right panel of \Cref{fig: double well sine perturbations} we can see the agreement of the sampled version with formula \eqref{eq: relative error double well time-dependent perturbation} when considering the time-dependent perturbation. We do not have a formula in the case of a space-dependent perturbation, however we can still observe the exponential dependence on the perturbation strength in the estimated relative error, which is expected for instance from formulas \eqref{eq: lower bound SDE KL} and \eqref{eqn: general bound relative error - 1}.

\subsection{Random stopping times}

\begin{figure}[H]
\centering
\includegraphics[width=0.40\linewidth]{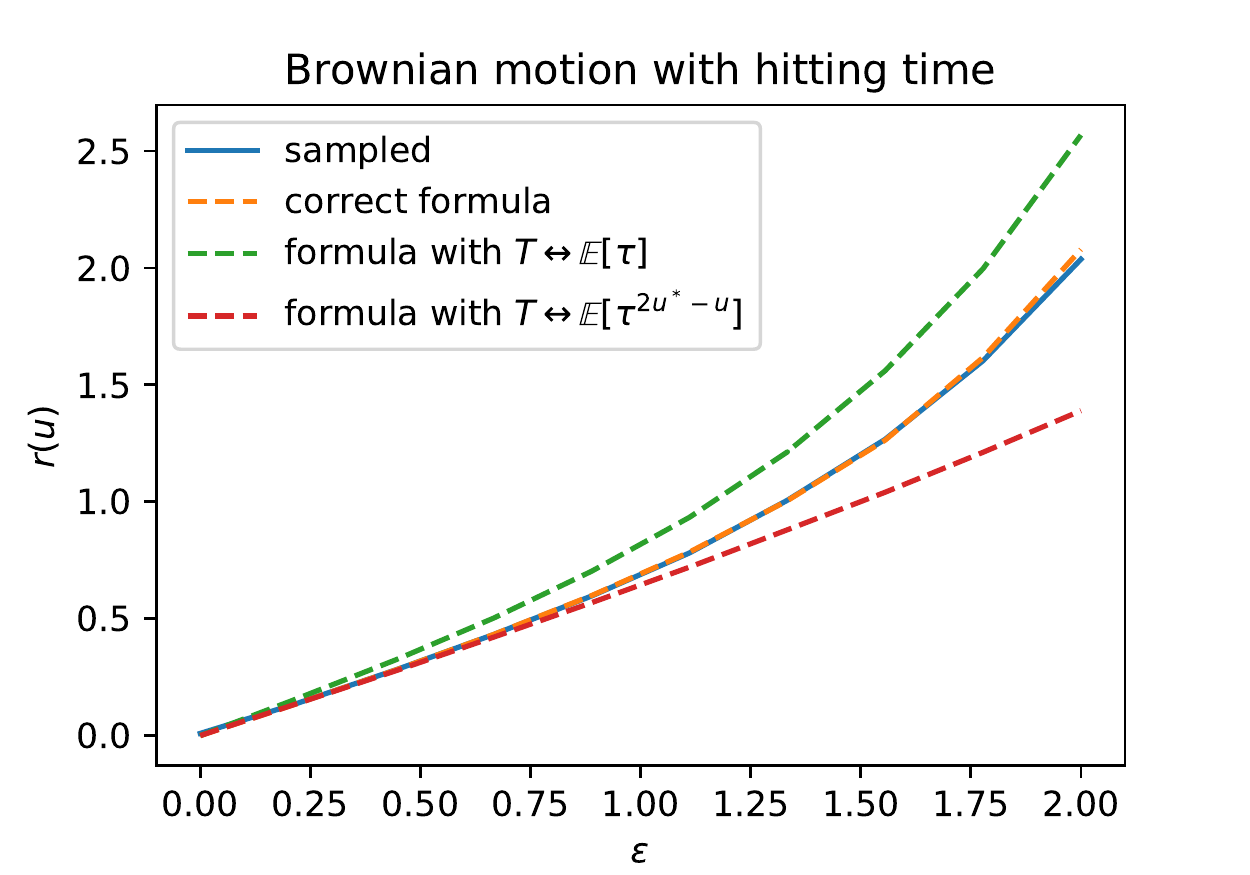}
\caption{Relative error of a quantity involving a random stopping time compared to the exact formula, a lower bound as well as a naive, but usually wrong approximation.}
\label{fig: brownian motion hitting 1d}
\end{figure}

The suboptimal importance sampling bounds from \Cref{sec: importance sampling in path space} can be transferred to problems that involve a random stopping time $\tau$ rather than a fixed time horizon $T$, where mostly $\tau^u = \inf \{ t > 0: X^u_t \notin \mathcal{D} \}$ is defined\footnote{We denote with $\tau = \tau^0$ the hitting time of the uncontrolled process $X_t$.} as the first exit time of a bounded domain $\mathcal{D} \subset \R^d$. However, one has to be careful with applying our formulas and bounds from above, as $\tau^u$ itself depends on the law of the process. For illustration, let us consider a one-dimensional toy example, where the dynamics is a scaled Brownian motion 
\begin{equation}
    X_t = \sqrt{2} W_t
\end{equation}
and we choose $f = 1, g = 0$ in \eqref{eqn: observable W}, such that  
\begin{equation}
    \mathcal{Z} = \E\left[e^{-\tau} \right].
\end{equation}
By noting that $\psi(x) = \E\left[e^{-\tau} | X_0 = x \right]$ fulfills the boundary value problem
\begin{subequations}
\begin{align}
    (\Delta - 1)\psi(x) &= 0, \qquad x \in \mathcal{D}, \\
    \psi(x) &= 1, \qquad x \in \partial \mathcal{D},
\end{align}
\end{subequations}
we can compute the optimal zero-variance importance sampling control to be
\begin{equation}
    u^*(x) = \sqrt{2} \nabla \log \psi(x) = \sqrt{2}\,\frac{1 - e^{-2x}}{e^{-2x} + 1}.
\end{equation}
In our experiment, we again perturb the optimal control via
\begin{equation}
    u = u^* + \varepsilon.
\end{equation}
Formula \eqref{eqn: general bound relative error - 2} provides an expression for the relative error, even if $T$ is replaced by a random time $\tau$ (which we leave the reader to check for herself), namely
\begin{equation}
\label{eqn: relative error hitting time}
    r(u) = \left(\E\left[e^{\varepsilon^2 \tau^{2u^* - u} } \right] - 1\right)^{\frac{1}{2}} \ge \left(e^{\varepsilon^2 \E\left[\tau^{2u^* - u} \right]}  - 1\right)^{\frac{1}{2}},
\end{equation}
where it is essential that $\tau^{2u^* - u}$ refers to the hitting time of the process $X_t^{2u^* - u}$. We applied Jensen's inequality in the last expression and note that naively assuming
\begin{equation}
\label{eqn: hitting time wrong formula}
     r(u) \approx \left(e^{\varepsilon^2 \E\left[\tau \right]}  - 1\right)^{\frac{1}{2}}
\end{equation}
is usually wrong. Figure \ref{fig: brownian motion hitting 1d} compares the sampled relative error with the exact formula, the lower bound in \eqref{eqn: relative error hitting time} and the wrong expression \eqref{eqn: hitting time wrong formula}.

\begin{remark}
Let us note again that estimating quantities involving hitting times gets particularly challenging in rare event settings, where the expected hitting time might become very large, cf. \Cref{ex: rare events of SDEs}. The relation \eqref{eqn: relative error hitting time} for the relative error then indicates that Monte Carlo estimation becomes especially difficult.
\end{remark}

\subsection{Small noise diffusions}
\label{sec: numerics small noise diffusion}

\begin{figure}[H]
\centering
\includegraphics[width=1.00\linewidth]{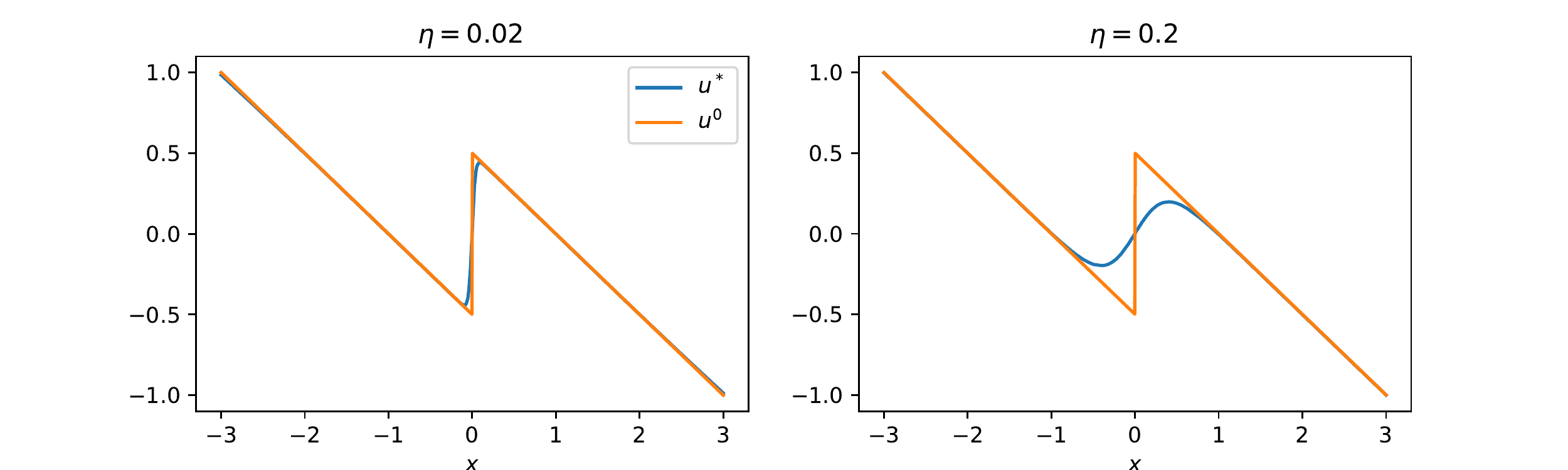}
\caption{For a small noise diffusion problem we display once the optimal control and once the control resulting from the zero-noise approximation with different noise scalings $\eta$.}
\label{fig: small noise diffusion V0 vs. V*}
\end{figure}

As an example for a small noise diffusion, we consider a modification of a one-dimensional toy example that has been proposed in \cite{vanden2012rare}. We take the scaled Brownian motion
\begin{equation}
    X_s^\eta = \sqrt{\eta} W_s, \qquad X_0 = 0.1,
\end{equation}
and want to compute 
\begin{equation}
    \E\left[e^{-\frac{1}{\eta}g(X_T^\eta)}\right]
\end{equation}
with 
\begin{equation}
    g(x) = \frac{\alpha}{2} \left(1-\frac{|x|}{\sqrt{\alpha}} \right)^2
\end{equation}
for $\alpha > 0$. One readily sees that
\begin{equation}
    V^0(x, t) = \frac{\alpha \left(1-\frac{|x|}{\sqrt{\alpha}} \right)^2}{2(T-t+1)}
\end{equation}
is the unique viscosity solution to the deterministic problem \eqref{eq: HJB deterministic problem}; we refer to \cite{fleming2006controlled} for a discussion of the theory of viscosity solutions. Since an explicit solution $V^*(x, t)$ to the second-order HJB equation \eqref{eq: HJB epsilon} is not available, we approximate it with finite differences.  In \Cref{fig: small noise diffusion V0 vs. V*} we show the corresponding controls $u^0(x, s) = -\sigma^\top V^0(x, t)$ and $u^*(x, s) = -\sigma^\top V^*(x, t)$ for different values of the noise coefficient $\eta$.

\begin{figure}[H]
\centering
\includegraphics[width=1.0\linewidth]{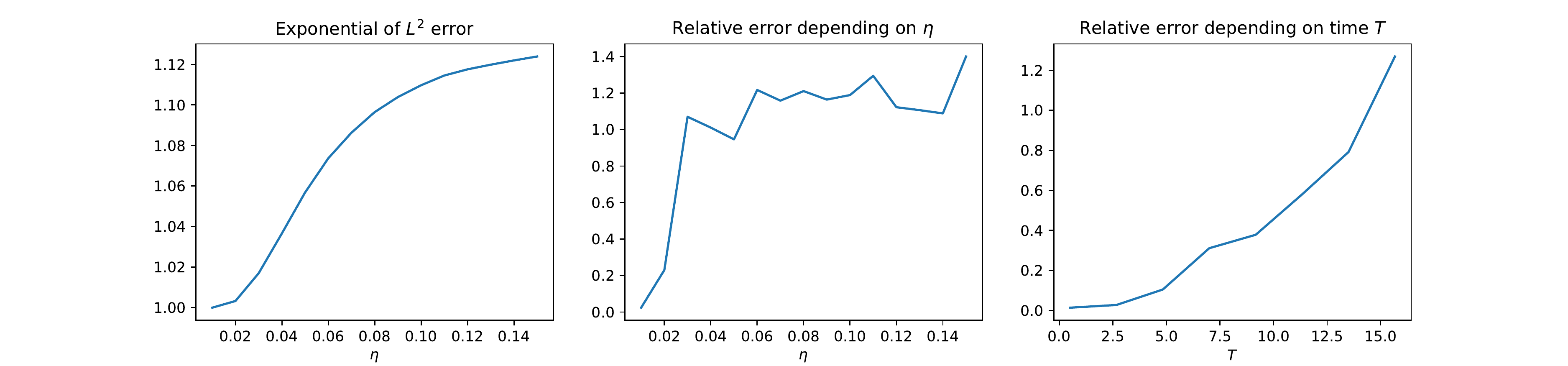}
\caption{Small noise diffusions with vanishing noise coefficient $\eta$. Left: Exponential of $L^2$ error between $u^*$ and $u^0$ depending on $\eta$ for $T = 1$. Middle: Relative importance sampling error depending on $\eta$. Right: Relative importance sampling error depending on $T$ for $\eta = 0.005$.}
\label{fig: small noise rel error T}
\end{figure}

In the middle panel of \Cref{fig: small noise rel error T} we show the relative error depending on the noise parameter $\eta$. Unlike one could expect from \eqref{eqn: relative error small noise diffusion}, it seems to not grow exponentially in $\eta$, which can be explained by looking at the exponentiated $L^2$ error, $\exp\left(\mathbb{E}\left[\int_0^T |u^* - u^0|^2(X^{u^0}_s, s) ds\right]\right)$, which we plot in the left panel. The observation that this does not grow exponentially seems to be rooted in the fact that the suboptimality $\delta = u^* - u^0$ is very different for different values of $x$. If we vary $T$, however, we can observe an exponential dependency on the time horizon, as displayed in the right panel of \Cref{fig: small noise rel error T}, again being in accordance with the consideration in \Cref{sec: small noise diffusions}.

\section{Conclusion and outlook}
\label{sec: conclusion}

In this article, we have provided quantitative bounds on the relative error of importance sampling that depend on the divergence between the actual proposal measure and the theoretically optimal one. These bounds indicate that importance sampling is very sensitive with respect to suboptimal choices of the proposals, which has been observed frequently in numerical experiments and is in line with recent theoretical analysis \cite{agapiou2015importance,sanz2018importance,chatterjee2018sample}. We showed that the relative error of importance sampling estimators scales exponentially in the KL divergence between the optimal and the proposal measure and argued that this renders importance sampling especially challenging in high dimensions.

We have focused on importance sampling of stochastic processes and derived some novel formulas for the relative error depending on the suboptimality of the function $u$ that controls the drift of the process. These formulas can be used to get practically useful bounds, but they also indicate two potential issues for importance sampling in path space: for systems with large state space and for problems on a long (or infinite) time horizon the relative error becomes exponentially large in the state space dimension $d$ and the time horizon $T$. We have briefly discussed how this observation can be transferred to random stopping times, such as first hitting times, and have applied our formulas to importance sampling in the small noise regime, offering new perspectives and revealing some potential drawbacks of existing methods.

Even though the key message of the paper regarding the use of path space importance sampling in high dimensions seems to be rather discouraging, let us finally mention that there is hope. In practice, the approximations to optimal proposals use iterative methods to minimize a divergence (or: loss function) between the approximant and the target, using, for instance, stochastic gradient descent. A crucial question then is which divergence to take, and it turns out that different choices lead to proposals with vastly different statistical properties. Let us mention four possible choices for the loss function in the approximation scheme: (a) relative entropy, (b) cross-entropy, (c) $\chi^2$ divergence or relative error, and (d) the recently introduced (see \cite{nusken2020solving}) log-variance divergence $\Var_{\widetilde{\nu}}(\log {\mathrm d}\nu^*/{\mathrm d}\widetilde{\nu})$, 

where we remark that in all four cases straightforward implementations for both probability densities on finite dimensional spaces and  (infinite-dimensional) path measures are available. Since normally we rely on Monte Carlo approximations of the measures and our quantities of interests, it is crucial that the relative error of these divergences and their gradients is as small as possible. While the analysis in this article suggests that the $\chi^2$ divergence %$\chi^2(\nu^*|\widetilde{\nu})$ 
cannot be expected to lead to a low-variance gradient estimator in general, the other divergences have been recently analyzed in \cite{nusken2020solving} in the context of path sampling (see \cite{richter2020vargrad} for related results on densities), some of which show better scaling properties when going to high dimensions. 

We expect that these perspectives can turn out fruitful in the future, in that they can guide the design of stable importance sampling schemes that work even in high dimensions. We therefore conclude that while importance sampling itself is often not robust, there are  strategies to approximate the optimal proposal measure in a more robust way that go beyond cross-entropy minimisation and control of the $\chi^2$ divergence.

\textbf{Acknowledgements}: This research has been funded by Deutsche Forschungsgemeinschaft (DFG) through the
grant CRC 1114 ‘Scaling Cascades in Complex Systems’ (project A05, project number 235221301). We
would like to thank  Wei Zhang and Nikolas N\"usken for many very useful discussions.

\newpage

\appendix
    
\section{Appendix}\label{app}

\subsection{Proofs for Section \ref{sec: importance sampling bounds divergences}}

\label{app: proofs - sec 2}

\begin{proof}[Proof of \Cref{prop: generalized Jensen}]
We adapt a proof of \cite{mitroi2011estimating}. Assume first that $m \ge 1$, then for any $E \in \mathcal{F}$
\begin{equation}
    \nu(E) - \lambda(E) \ge \nu(E) - m\lambda(E) \ge 0,
\end{equation}
where the last inequality follows from the definition of $m$. On the other hand, if $E = \Omega$, then $\nu(E) - \lambda(E) = 0$ and therefore it follows that $m = 1$, i.e. $\nu = \lambda$.\par\bigskip
Let now $m < 1$. We want to show $m J(f, \lambda, \varphi) \le J(f, \nu, \varphi)$, which is equivalent to
\begin{equation}
     \E_{\nu}[f(\varphi)] - m\E_{\lambda}[f(\varphi)] + m f\left(\E_{\lambda}[\varphi]\right) \ge f\left(\E_{\nu}[\varphi] \right).
\end{equation}
We compute
\begin{subequations}
\begin{align}
    \E_{\nu}[f(\varphi)]  - m\E_{\lambda}[f(\varphi)] + m f\left(\E_{\lambda}[f(\varphi)] \right) &\ge \left(\E_{\nu}[1] - m\E_{\lambda}[1] \right) f\left(\frac{\E_{\nu}[\varphi] -m\E_{\lambda}[\varphi]}{\E_{\nu}[1] - m\E_{\lambda}[1]}\right)+ m f\left(\E_{\lambda}[f(\varphi)] \right) \\
    &= (1-m)f\left(\frac{\E_{\nu}[\varphi] -m\E_{\lambda}[\varphi]}{1-m}  \right) + m f\left(\E_{\lambda}[f(\varphi)] \right) \\
    &\ge f\left(\E_{\nu}[\varphi]  - m\E_{\lambda}[\varphi] + m \E_{\lambda}[\varphi] \right) \\
    &= f\left(\E_{\nu}[\varphi] \right),
\end{align}
\end{subequations}
where we used two times the convexity of $f$. The other inequality follows analogously.
\end{proof}

\subsection{Proofs for Section \ref{sec: importance sampling in path space}}

\label{app: proofs}

\begin{proof}[Proof of Proposition \ref{prop: general bound relative error}]
We compute
\begin{subequations}
\begin{align}
    \E\left[ e^{-2 \W(X^u)} \left( \frac{\mathrm d \P}{\mathrm d \P^u}(X^u) \right)^2 \right] &= \E\left[ e^{-2 \W(X^u)} \left( \frac{\mathrm d \P}{\mathrm d \P^{u^*}}(X^u) \frac{\mathrm d \P^{u^*}}{\mathrm d \P^u}(X^u) \right)^2 \right] \\
    &= {\Z}^2 \E \left[ \left(\frac{\mathrm d \P^{u^*}}{\mathrm d \P^u}(X^u) \right)^2 \right],
\end{align}
\end{subequations}
where we used
\begin{equation}
    \frac{\mathrm d \P}{\mathrm d \P^{u^*}}(X^u) = e^{\W(X^u)} \Z.
\end{equation}
Equation \eqref{eqn: general bound relative error - 1} now follows by the Girsanov formula (see \Cref{lem: Girsanov}) and the definition of the variance. For equation \eqref{eqn: general bound relative error - 2} note that we can write
\begin{align}
    {\E}_{\P^u} \left[ \left(\frac{\mathrm d \P^{u^*}}{\mathrm d \P^u} \right)^2 \right] =  {\E}_{\P^{u + 2 \delta}} \left[ \left(\frac{\mathrm d \P^{u^*}}{\mathrm d \P^u} \right)^2 \frac{\mathrm d \P^{u}}{\mathrm d \P^{u + 2\delta}} \right].
\end{align}
We compute 
\begin{align}
    \frac{\mathrm d \P^{u^*}}{\mathrm d \P^u}(X^{u+2\delta}) = \exp\left( \frac{3}{2} \int_0^T |\delta(X_s^{u+2\delta}, s)|^2\mathrm ds + \int_0^T \delta(X_s^{u + 2\delta}, s) \cdot \mathrm d W_s \right)
\end{align}
and
\begin{align}
    \frac{\mathrm d \P^{u}}{\mathrm d \P^{u+2\delta}}(X^{u+2\delta}) = \exp\left( -2 \int_0^T |\delta(X_s^{u+2\delta}, s)|^2\mathrm ds - 2\int_0^T \delta(X_s^{u + 2\delta}, s) \cdot \mathrm d W_s \right),
\end{align}
 from which the desired formula immediately follows.
\end{proof}

\begin{proof}[Alternative proof of \Cref{prop: RE bounds independent of x}]

We follow the reasoning in \cite[Thm.~2.1]{klebaner2014benes} and apply Gr\"onwall's inequality to the square integrable exponential martingale $Z$.\footnote{See also Theorem 2 in \url{http://math.ucsd.edu/~pfitz/downloads/courses/spring05/math280c/expmart.pdf}.} To this end, we define the shorthands $\delta(x, t) := (u^* - u)(x, t)$ and
\begin{equation}
    Z_t := \exp\left(-\frac{1}{2} \int_0^t |\delta(X_s, s)|^2 \mathrm ds + \int_0^t \delta(X_s, s)\cdot\mathrm dW_s\right).
\end{equation}
Then, by It\^{o}'s formula, 
\begin{equation}
    Z_t^2 = 1 + 2\int_0^t Z_s \mathrm dZ_s + \int_0^t Z_s^2 |\delta(X_s, s)|^2 \mathrm ds\,
\end{equation}
and therefore, after taking expectations,  
\begin{subequations}
\begin{align}
    \E\left[Z^2_t \right] &= 1 + \E\left[ \int_0^t Z_s^2 |\delta(X_s, s)|^2 \mathrm ds \right] \\
    &\le 1 +  \int_0^t \E\left[Z_s^2\right] h_2^2(s) \mathrm ds.
\end{align}
\end{subequations}
We can now apply Gr\"onwall's inequality to get
\begin{equation}
    \E\left[Z^2_t \right] \le \exp\left(\int_0^t h_2^2(s)\mathrm ds \right)
\end{equation}
and therefore the desired statement after applying Proposition \ref{prop: general bound relative error}. The other direction follows analogously by noting that 
\begin{align}
    -\E\left[Z^2_t \right] \le -1 -  \int_0^t \E\left[Z_s^2\right] h_1^2(s) \mathrm ds.
\end{align}
\end{proof}
\begin{remark} Yet another alternative to prove \Cref{prop: RE bounds independent of x} is by computing 

\begin{subequations}
\begin{align}
{\E}\left[\left(\frac{\mathrm d \P^{u^*}}{\mathrm d \P^u}(X^u) \right)^2 \right] &= {\E}\left[\exp\left(-\int_0^T |\delta(X^u_s, s)|^2 \mathrm ds + 2 \int_0^T \delta(X^u_s, s) \cdot \mathrm dW_s \right)\right] \\ 
&= {\E}\left[\exp\left(\int_0^T |\delta(X^u_s, s)|^2 \mathrm ds - 2 \int_0^T |\delta(X^u_s, s)|^2 \mathrm ds + 2 \int_0^T \delta(X^u_s, s) \cdot \mathrm dW_s \right)\right] \\ 
&\le \exp\left( \int_0^T h_2^2(s) \mathrm ds \right)\, {\E}\left[\exp\left(- \frac{1}{2} \int_0^T | 2\delta(X^u_s, s)|^2 \mathrm ds + \int_0^T 2\delta(X^u_s, s) \cdot \mathrm dW_s \right)\right] \\
&= \exp\left( \int_0^T h_2^2(s) \mathrm ds \right),
\end{align}
\end{subequations}
where we used the constant expectation property of the exponential martingale in the last step. The other direction follows analogously.
\end{remark}

\begin{proof}[Proof of Proposition \ref{prop: Hoelder bound}]
\label{proof: Hoelder bound}
From Lemma \ref{lemma: bound with Hoelder} it holds for $n, p, q > 1$ with $\frac{1}{p} + \frac{1}{q} = 1$ that
\begin{equation}
    {\E}\left[\left(\frac{\mathrm d \P^{u^*}}{\mathrm d \P^u}(X^u) \right)^n \right] \le {\E}\left[\exp \left(  \frac{nq(np-1)}{2} \int_0^T |u^* - u|^2(X^u_s,s) \mathrm ds  \right) \right]^{\frac{1}{q}}.
\end{equation}
We write $q = \frac{p}{p-1}$ and note that $q(np-1) = \frac{p(np-1)}{p-1}$ is minimized by $p^* = 1 \pm \sqrt{1 - \frac{1}{n}}$, from which we are only allowed to take the positive part due to the constraint $p \ge 1$. For $n=2$ this yields $p^* = \frac{\sqrt{2}+1}{\sqrt{2}}$ and $q^*=\sqrt{2}+1$, and we get the desired statement by recalling

\begin{equation}
    r^2(u) = {\Var}_{\P^u}\left( \frac{\mathrm d \P^{u^*}}{\mathrm d \P^u}\right) = {\E}_{{\P}^{u}}\left[\left(\frac{\mathrm d \P^{u^*}}{\mathrm d \P^u} \right)^2 \right] - 1.
\end{equation}
\end{proof}

\subsection{Auxiliary statements}

In this section, we recall some known statements and provide some helpful additional analysis.\par\bigskip

First note that the Radon-Nikodym derivative appearing in the importance sampling estimator in path space can be computed explicitly.
\begin{lemma}[Girsanov]
\label{lem: Girsanov}
For $u \in \mathcal{U}$, the measures $\mathbb{P}$ and $\mathbb{P}^u$, relating to the SDEs \eqref{eq: uncontrolled SDE} and \eqref{eq: controlled SDE}, are equivalent. Moreover, the Radon-Nikodym derivative satisfies
\begin{equation}
\label{eq:Pu P}
\frac{\mathrm d \P^{u}}{\mathrm d\P}(X) = \exp \left( \int_0^T \left(u^\top  \sigma^{-1}\right)(X_s, s) \cdot \mathrm dX_s - \int_0^T (\sigma^{-1} b \cdot u)(X_s, s) \,\mathrm ds - \frac{1}{2} \int_0^T |u(X_s, s)|^2 \, \mathrm ds \right).
\end{equation}
\end{lemma}

\begin{proof}
See \cite[Lemma A.1]{nusken2020solving}.
\end{proof}

\begin{corollary}[Formula for path space relative error in a special case]
If the difference $u^* - u$ does not depend on $x$, then
\begin{equation}
    r(u) = \left(\exp\left(\int_0^T |u^* - u|^2(s) \mathrm ds\right) - 1\right)^{\frac{1}{2}}.
\end{equation}
\end{corollary}

\begin{proof}
This is a direct consequence of \eqref{eqn: general bound relative error - 2}. For the reader's convenience, we provide an alternative proof. If $u^* - u$ does not depend on $x$, then the random variable 
\begin{equation}
    Y = -\int_0^T |u^* - u|^2(s) \mathrm ds + 2 \int_0^T (u^* - u)(s) \cdot \mathrm dW_s 
\end{equation}
is normally distributed, with mean and variance given by 
\begin{equation}
    \mu = -\int_0^T |u^* - u|^2(s) \mathrm ds, \qquad \sigma^2 = 4\int_0^T |u^* - u|^2(s) \mathrm ds,
\end{equation}
where the second expression follows from the It\^{o} isometry. The random variable $\left(\frac{\mathrm d \P^{u^*}}{\mathrm d \P^u}(X^u) \right)^2 = e^Y$ is then log-normally distributed and we compute
\begin{equation}
    \E\left[ e^Y \right] = e^{\mu + \frac{\sigma^2}{2}} = e^{\frac{\sigma^2}{4}},
\end{equation}
which gives the desired statement.
\end{proof}

\begin{lemma}
\label{lemma: bound with Hoelder}
Let $n, p, q > 1$ with $\frac{1}{p} + \frac{1}{q} = 1$, then it holds that
\begin{equation}
    {\E}\left[\left(\frac{\mathrm d \P^{u^*}}{\mathrm d \P^u}(X^u) \right)^n \right] \le {\E}\left[\exp \left(  \frac{nq(np-1)}{2} \int_0^T |u^* - u|^2(X^u_s,s) \mathrm ds  \right) \right]^{\frac{1}{q}}.
\end{equation}
\end{lemma}
\begin{proof}

Let us write $\delta(x, s) := (u^*-u)(x, s)$, and let $n, p, q > 1$, then, using the H\"older inequality with $\frac{1}{p} + \frac{1}{q} = 1$, it holds
\begin{subequations}
\begin{align}
{\E}_{{\P}^{u}}\left[\left(\frac{\mathrm d \P^{u^*}}{\mathrm d \P^u} \right)^n \right] &= {\E}_{\P^u}\left[ \exp\left(n\int_0^T \delta(X_s, s) \cdot \mathrm d W_s - \frac{n^2p}{2} \int_0^T |\delta(X_s,s)|^2 \mathrm ds  + \frac{n(np-1)}{2} \int_0^T |\delta(X_s,s)|^2 \mathrm ds  \right) \right] \\
&\le {\E}_{\P^u}\left[ \exp\left(\int_0^T np \, \delta(X_s, s) \cdot \mathrm d W_s - \frac{1}{2}\int_0^T |np \, \delta(X_s,s)|^2 \mathrm ds \right) \right]^{\frac{1}{p}} \\
& \,\,\quad{\E}_{\P^u}\left[\exp \left(  \frac{nq(np-1)}{2} \int_0^T |\delta(X_s,s)|^2 \mathrm ds  \right) \right]^{\frac{1}{q}} \\
&= {\E}_{\P^u}\left[\exp \left(  \frac{nq(np-1)}{2} \int_0^T |\delta(X_s,s)|^2 \mathrm ds  \right) \right]^{\frac{1}{q}}.
\end{align}
\end{subequations}
Note that, even though H\"older's inequality holds for $p, q \in [1, \infty]$, the inequality becomes useless for $q = 1$ and $p = \infty$.
\end{proof}

\begin{proposition}[Zero-variance property]
\label{prop: u^* implies zero-variance}
We get a vanishing relative error $r(u) = 0$ if and only if $\delta = u^* - u = 0$, i.e. when having the optimal control $u = u^* = -\sigma^\top \nabla V$.
\end{proposition}
\begin{proof}
The fact that $\delta = 0$ implies $r(u) = 0$ follows directly from \eqref{eqn: general bound relative error - 2} or \eqref{eqn: moment term h(x,t)}. For the other direction note that $r(u) = 0$ implies $M_{u}(x, t) = \psi^2(x, t)$ (as defined in \Cref{prop: PDE for moment ratio}) for all $(x, t) \in \R^d \times [0, T]$ and therefore equation \eqref{eqn: PDE second moment} becomes
\begin{equation}
    (\partial_t + L - \sigma u(x, t) \cdot \nabla - 2 f(x, t) + |u(x, t)|^2) \psi^2(x, t) = 0.
\end{equation}
Further note that due to the Kolmogorov backward equation it holds
\begin{equation}
    (\partial_t + L - 2f(x, t))\psi^2(x, t) - |(\sigma^\top\nabla \psi)(x, t)|^2 = 0.
\end{equation}
Combining these two PDEs brings
\begin{equation}
     \psi^2(x, t)|u(x, t)|^2 - 2 (\psi \sigma u \cdot \nabla \psi)(x, t) + |(\sigma^\top \nabla \psi)(x, t)|^2 = |(\psi u)(x, t) - (\sigma^\top \nabla \psi)(x, t)|^2 = 0,
\end{equation}
which implies that 
\begin{equation}
    u = \sigma^\top \frac{\nabla \psi}{\psi} = \sigma^\top  \nabla \log \psi=- \sigma^\top  \nabla V.
\end{equation}
\end{proof}

The following lemma shows that the $\mathrm{KL}$ divergence increases with the number of dimensions.  This result follows from the chain-rule of KL divergence, see, e.g., \cite{cover2012elements}.

\begin{lemma}[Dimension dependence of KL divergence]
\label{lem: KL dimension dependence}
Let $u^{(d)}(z_1,\ldots,z_d)$ and $v^{(d)}(z_1,\ldots,z_d)$ be two arbitrary probability distributions on $\mathbb{R}^d$. For $j \in  \{1 \ldots,d \}$ denote their marginals on the first $j$ coordinates by $u^{(j)}$ and $v^{(j)}$, i.e. 
\begin{equation}
u^{(j)}(z_1,\ldots,z_j) = \idotsint u^{(d)}(z_1, \ldots, z_d) \, \mathrm{d}z_{j+1} \ldots \mathrm{d}z_d, 
\end{equation}
and 
\begin{equation}
v^{(j)}(z_1,\ldots,z_j) = \idotsint v^{(d)}(z_1, \ldots, z_d) \, \mathrm{d}z_{j+1} \ldots \mathrm{d}z_d.
\end{equation}
Then \begin{equation}
\mathrm{KL}(u^{(1)} \; | \; v^{(1)}) \le \mathrm{KL}(u^{(2)} \; | \; v^{(2)}) \le \ldots \le \mathrm{KL}(u^{(d)} \; | \; v^{(d)}),
\end{equation}
i.e. the function $J \mapsto \mathrm{KL}(u^{(j)} \; | \; v^{(j)})$ is increasing.
\end{lemma}

\subsubsection{Relative error of log-normal random variables}
\label{app: relative error log-normal}
Let $Y \sim \mathcal{N}(\overline{\mu}, \overline{\Sigma})$ with arbitrary $\overline{\mu} \in \R^d, \overline{\Sigma} \in \R^{d\times d}$ and take $\gamma \in \R^d, c\in \R$, then $e^{\gamma \cdot Y + c}$ is log-normally distributed and its relative error is
\begin{equation}
    r(\gamma, \overline{\Sigma}) = \sqrt{\frac{\E\left[e^{2 (\gamma \cdot Y + c)} \right]}{\E\left[e^{\gamma \cdot Y + c} \right]^2} - 1}= \sqrt{\frac{\E\left[e^{2 \gamma \cdot Y} \right]}{\E\left[e^{\gamma \cdot Y } \right]^2} - 1} =  \sqrt{e^{\gamma\cdot \overline{\Sigma} \gamma}-1},
\end{equation}
independent of $c$.
With the setting and notation from \Cref{ex: high-dim Gaussians} we can now for instance compute
\begin{equation}
e^{-g(\widetilde{X})}\frac{p}{\widetilde{p}^\varepsilon}(\widetilde{X}) =\exp\left( - \alpha \cdot \widetilde{X} + \log\frac{p}{\widetilde{p}^\varepsilon}\right) =\exp\left( \varepsilon \cdot \widetilde{X} - \mu\cdot(\alpha + \varepsilon) + \frac{1}{2}(\alpha + \varepsilon)\cdot \Sigma (\alpha + \varepsilon)\right)
\end{equation} 
and with $\gamma = \varepsilon, c = - \mu\cdot(\alpha + \varepsilon) +\frac{1}{2}(\alpha + \varepsilon)\cdot \Sigma (\alpha + \varepsilon), \overline{\Sigma} = \Sigma$ one therefore gets the relative error
\begin{equation}
        r(\widetilde{p}^\varepsilon) = \sqrt{e^{\varepsilon\cdot {\Sigma} \varepsilon}-1}
\end{equation}
as stated in \eqref{eq: r for epsilon shifted Gaussian}.

\subsubsection{Asymptotic expansion in small noise diffusions}
\label{app: asymptotic expansion small noise}
To get further intuition on the small noise diffusions defined in \Cref{sec: small noise diffusions}, let us consider the formal expansion of the solution to the HJB equation \eqref{eq: HJB epsilon}
\begin{equation}
    V = v_0 + \eta v_1 + \eta^2 v_2 + \dots.
\end{equation}
Inserting into \eqref{eq: HJB epsilon} (with $\sigma = I_{d \times d}$) and comparing the powers of $\eta$ yields the PDEs
\begin{subequations}
\begin{align}
    &\partial_t v_0 + b \cdot \nabla v_0 - \frac{1}{2}\left| \nabla v_0 \right|^2 = 0, \\
    &\partial_t v_1 + \frac{1}{2} \Delta v_0 + b \cdot \nabla v_1 -  \nabla v_0 \cdot \nabla v_1 = 0, \\
    &\partial_t v_2 + \frac{1}{2} \Delta v_1 + b \cdot \nabla v_2 -  \nabla v_0 \cdot \nabla v_2 - \frac{1}{2} \left|\nabla v_1\right|^2 = 0,
\end{align}
\end{subequations}
and so on, where all but the first PDE are transport equations (see \cite{spiliopoulos2015nonasymptotic}). We note that (given some appropriate assumptions) we have $v_0 = V^0$, with $V^0$ being to solution to \eqref{eq: HJB deterministic problem}. In \cite{fleming1971stochastic} it is proven that 
\begin{equation}
     \nabla V = \nabla  V^0 + \eta \nabla  v_1 + o(\eta),
\end{equation} 
where $v_1$ fulfills the PDE above and $V$ is the solution to the original HJB equation \eqref{eq: HJB epsilon}.

\printbibliography{}

@incollection{bengtsson2008curse,
  title={Curse-of-dimensionality revisited: Collapse of the particle filter in very large scale systems},
  author={Bengtsson, Thomas and Bickel, Peter and Li, Bo and others},
  booktitle={Probability and statistics: Essays in honor of David A. Freedman},
  pages={316--334},
  year={2008},
  publisher={Institute of Mathematical Statistics}
}

@book{dembo2009large,
  title={Large Deviations Techniques and Applications},
  author={Dembo, Amir and Zeitouni, Ofer},
  year={2009},
  publisher={Springer Berlin Heidelberg}
}

@article{bengtsson2005curse,
	title={Curse-of-dimensionality revisited: {C}ollapse of importance sampling in very high-dimensional systems},
	author={Li, Bo and Bengtsson, Thomas and Bickel, Peter},
	pages={1-18},
	volume = {696},
	year={2005},
	journal={Tech Reports, Department of Statistics, UC Berkeley}
}

@article{klebaner2014benes,
author = {Klebaner, Fima and Liptser, Robert},
title = {When a Stochastic Exponential Is a True Martingale. {E}xtension of the {B}ene\v{s} Method},
journal = {Theory of Probability and its Applications},
volume = {58},
number = {1},
pages = {38-62},
year = {2014}
}

@article{Bovier2005,
author = {Bovier, Anton and Gayrard, V\'eronique and Klein, Markus},
journal = {J. Eur. Math. Soc.},
number = {1},
volume={7}, 
pages = {69-99}, 
title = {Metastability in reversible diffusion processes {II}: precise asymptotics for small eigenvalues},
year = {2005}
}

@article{chatterjee2018sample,
  title={The sample size required in importance sampling},
  author={Chatterjee, Sourav and Diaconis, Persi},
  journal={The Annals of Applied Probability},
  volume={28},
  number={2},
  pages={1099--1135},
  year={2018},
  publisher={Institute of Mathematical Statistics}
}

@article{sason2014improved,
  title={On Improved Bounds for Probability Metrics and $f$-Divergences},
  author={Sason, Igal},
  journal={arXiv preprint arXiv:1403.7164},
  year={2014}
}

@article{dragomir2000some,
  title={Some Inequalities For The {K}ullback-{L}eibler And $\chi^2$-Distances In Information Theory And Applications},
  author={Dragomir, Sever S and Gluscevic, V},
  journal={RGMIA research report collection},
  volume={3},
  number={2},
  pages={199--210},
  year={2000},
  publisher={School of Communications and Informatics, Faculty of Engineering and Science~…}
}

@article{spiliopoulos2015nonasymptotic,
  title={Nonasymptotic performance analysis of importance sampling schemes for small noise diffusions},
  author={Spiliopoulos, Konstantinos},
  journal={Journal of Applied Probability},
  volume={52},
  number={3},
  pages={797--810},
  year={2015},
  publisher={Cambridge University Press}
}

@article{vanden2012rare,
  title={Rare event simulation of small noise diffusions},
  author={Vanden-Eijnden, Eric and Weare, Jonathan},
  journal={Communications on Pure and Applied Mathematics},
  volume={65},
  number={12},
  pages={1770--1803},
  year={2012},
  publisher={Wiley Online Library}
}

@article {Ragone2018,
	author = {Ragone, Francesco and Wouters, Jeroen and Bouchet, Freddy},
	title = {Computation of extreme heat waves in climate models using a large deviation algorithm},
	volume = {115},
	number = {1},
	pages = {24-29},
	year = {2018},
	doi = {10.1073/pnas.1712645115},
	journal = {Proceedings of the National Academy of Sciences}
}

@article{glasserman1997counterexamples,
  title={Counterexamples in importance sampling for large deviations probabilities},
  author={Glasserman, Paul and Wang, Yashan},
  journal={The Annals of Applied Probability},
  volume={7},
  number={3},
  pages={731--746},
  year={1997},
  publisher={Institute of Mathematical Statistics}
}

@book{asmussen2007stochastic,
  title={Stochastic Simulation: Algorithms and Analysis},
  author={Asmussen, S\o ren and Glynn, Peter W.},
  year={2007},
  publisher={Springer, New York}
}

@article{asmussen2011importance,
  title={Importance sampling for rare events},
  author={Asmussen, S\o ren and Dupuis, Paul and Rubinstein, Reuven and Wang, Hui},
  journal={Aarhus Univ., Aarhus, Denmark, Tech. Rep.},
  year={2011},
  publisher={Citeseer}
}

@article{dupuis2015escaping,
  title={Escaping from an attractor: Importance sampling and rest points {I}},
  author={Dupuis, Paul and Spiliopoulos, Konstantinos and Zhou, Xiang},
  journal={The Annals of Applied Probability},
  pages={2909--2958},
  year={2015},
  publisher={JSTOR}
}

@article{nusken2020solving,
  title={Solving high-dimensional {H}amilton-{J}acobi-{B}ellman {PDEs} using neural networks: perspectives from the theory of controlled diffusions and measures on path space},
  author={N{\"u}sken, Nikolas and Richter, Lorenz},
  journal={arXiv preprint arXiv:2005.05409},
  year={2020}
}

@article{fleming1971stochastic,
  title={Stochastic control for small noise intensities},
  author={Fleming, Wendell H.},
  journal={SIAM Journal on Control},
  volume={9},
  number={3},
  pages={473--517},
  year={1971},
  publisher={SIAM}
}

@article{hartmann2017variational,
  title={Variational characterization of free energy: Theory and algorithms},
  author={Hartmann, Carsten and Richter, Lorenz and Sch{\"u}tte, Christof and Zhang, Wei},
  journal={Entropy},
  volume={19},
  number={11},
  pages={626},
  year={2017},
  publisher={Multidisciplinary Digital Publishing Institute}
}

@article{mitroi2011estimating,
  title={Estimating the normalized {J}ensen functional},
  author={Mitroi, Flavia Corina},
  journal={J. Math. Inequal},
  volume={5},
  number={4},
  pages={507--521},
  year={2011}
}

@inproceedings{sason2015tight,
  title={Tight bounds for symmetric divergence measures and a new inequality relating $f$-divergences},
  author={Sason, Igal},
  booktitle={2015 IEEE Information Theory Workshop (ITW)},
  pages={1--5},
  year={2015},
  organization={IEEE}
}

@book{cover2012elements,
  title={Elements of {I}nformation {T}heory},
  author={Cover, Thomas M and Thomas, Joy A},
  year={2012},
  publisher={John Wiley \& Sons}
}

@article{deboer2005tutorial,
  title={A tutorial on the cross-entropy method},
  author={De Boer, Pieter-Tjerk and Kroese, Dirk P and Mannor, Shie and Rubinstein, Reuven Y},
  journal={Annals of operations research},
  volume={134},
  number={1},
  pages={19--67},
  year={2005},
  publisher={Springer}
}

@article{zhang2014applications,
  title={Applications of the cross-entropy method to importance sampling and optimal control of diffusions},
  author={Zhang, Wei and Wang, Han and Hartmann, Carsten and Weber, Marcus and Sch\"utte, Christof},
  journal={SIAM Journal on Scientific Computing},
  volume={36},
  number={6},
  pages={A2654--A2672},
  year={2014},
  publisher={SIAM}
}

@article{berglund2011kramers,
  title={Kramers' law: Validity, derivations and generalisations},
  author={Berglund, Nils},
  journal={Markov Processes and Related fields},
  volume = {19},
  number = {3},
  pages = {459-490},
  year={2013}
}

@article{richter2020vargrad,
  title={Var{G}rad: A Low-Variance Gradient Estimator for Variational Inference},
  author={Richter, Lorenz and Boustati, Ayman and N{\"u}sken, Nikolas and Ruiz, Francisco JR and Akyildiz, {\"O}mer Deniz},
  journal={arXiv preprint arXiv:2010.10436},
  year={2020}
}

@article{meng1996simulating,
  title={Simulating ratios of normalizing constants via a simple identity: a theoretical exploration},
  author={Meng, Xiao-Li and Wong, Wing Hung},
  journal={Statistica Sinica},
  pages={831--860},
  year={1996},
  publisher={JSTOR}
}

@book{fleming2006controlled,
  title={Controlled Markov processes and viscosity solutions},
  author={Fleming, Wendell H. and Soner, Halil Mete},
  volume={25},
  year={2006},
  publisher={Springer Science \& Business Media}
}

@article{spiliopoulos2013large,
  title={Large deviations and importance sampling for systems of slow-fast motion},
  author={Spiliopoulos, Konstantinos},
  journal={Applied Mathematics \& Optimization},
  volume={67},
  number={1},
  pages={123--161},
  year={2013},
  publisher={Springer}
}

@article{dupuis2012importance,
  title={Importance sampling for multiscale diffusions},
  author={Dupuis, Paul and Spiliopoulos, Konstantinos and Wang, Hui},
  journal={Multiscale Modeling \& Simulation},
  volume={10},
  number={1},
  pages={1--27},
  year={2012},
  publisher={SIAM}
}

@book{owen2013monte,
  title={Monte Carlo theory, methods and examples},
  author={Owen, Art B.},
  publisher={Self-published},
  year={2013}
}

@book{glasserman2013monte,
  title={Monte Carlo methods in financial engineering},
  author={Glasserman, Paul},
  volume={53},
  year={2013},
  publisher={Springer Science \& Business Media}
}

@book{delmoral2013ips,
  title={Mean Field Simulation for {M}onte {C}arlo Integration},
  author={Del Moral, Pierre},
  year={2013},
  publisher={Chapman and Hall/CRC}
}

@article{dupuis2004importance,
  title={Importance sampling, large deviations, and differential games},
  author={Dupuis, Paul and Wang, Hui},
  journal={Stochastics: An International Journal of Probability and Stochastic Processes},
  volume={76},
  number={6},
  pages={481--508},
  year={2004},
  publisher={Taylor \& Francis}
}

@article{dupuis2007subsolutions,
  title={Subsolutions of an {I}saacs equation and efficient schemes for importance sampling},
  author={Dupuis, Paul and Wang, Hui},
  journal={Mathematics of Operations Research},
  volume={32},
  number={3},
  pages={723--757},
  year={2007},
  publisher={INFORMS}
}

@book{liu2008monte,
  title={Monte Carlo strategies in scientific computing},
  author={Liu, Jun S},
  year={2008},
  publisher={Springer Science \& Business Media}
}

@article{siegmund1976importance,
  title={Importance sampling in the {M}onte {C}arlo study of sequential tests},
  author={Siegmund, David},
  journal={The Annals of Statistics},
  pages={673--684},
  year={1976},
  publisher={JSTOR}
}

@incollection{doucet2001introduction,
  title={An introduction to sequential {M}onte {C}arlo methods},
  author={Doucet, Arnaud and De Freitas, Nando and Gordon, Neil},
  booktitle={Sequential {M}onte {C}arlo methods in practice},
  pages={3--14},
  year={2001},
  publisher={Springer}
}

@article{hartmann2012efficient,
  title={Efficient rare event simulation by optimal nonequilibrium forcing},
  author={Hartmann, Carsten and Sch{\"u}tte, Christof},
  journal={Journal of Statistical Mechanics: Theory and Experiment},
  volume={2012},
  number={11},
  pages={P11004},
  year={2012},
  publisher={IOP Publishing}
}

@article{hartmann2019chaos,
author = {Hartmann,Carsten  and Kebiri,Omar  and Neureither,Lara  and Richter,Lorenz },
title = {Variational approach to rare event simulation using least-squares regression},
journal = {Chaos},
volume = {29},
number = {6},
pages = {063107},
year = {2019},
doi = {10.1063/1.5090271}
}

@book{bishop2006pattern,
  title={Pattern recognition and machine learning},
  author={Bishop, Christopher M},
  year={2006},
  publisher={Springer}
}

@book{stoltz2010free,
  title={Free energy computations: A mathematical perspective},
  author={Stoltz, Gabriel and Rousset, Mathias and others},
  year={2010},
  publisher={World Scientific}
}

@article{gelman1998simulating,
  title={Simulating normalizing constants: From importance sampling to bridge sampling to path sampling},
  author={Gelman, Andrew and Meng, Xiao-Li},
  journal={Statistical science},
  pages={163--185},
  year={1998},
  publisher={JSTOR}
}

@article{glasserman2005importance,
  title={Importance sampling for portfolio credit risk},
  author={Glasserman, Paul and Li, Jingyi},
  journal={Management science},
  volume={51},
  number={11},
  pages={1643--1656},
  year={2005},
  publisher={INFORMS}
}

@article{agapiou2015importance,
  title={Importance sampling: computational complexity and intrinsic dimension},
  author={Agapiou, Sergios and Papaspiliopoulos, Omiros and Sanz-Alonso, Daniel and Stuart, Andrew M.},
  journal={Statistical Science},
  volume = {32},
  number = {3},
  year={2015}
}

@article{sanz2018importance,
  title={Importance sampling and necessary sample size: an information theory approach},
  author={Sanz-Alonso, Daniel},
  journal={SIAM/ASA Journal on Uncertainty Quantification},
  volume={6},
  number={2},
  pages={867--879},
  year={2018},
  publisher={SIAM}
}

@article{polyak2017does,
  title={Why does {M}onte {C}arlo fail to work properly in high-dimensional optimization problems?},
  author={Polyak, Boris and Shcherbakov, Pavel},
  journal={Journal of Optimization Theory and Applications},
  volume={173},
  number={2},
  pages={612--627},
  year={2017},
  publisher={Springer}
}

@article{chen2005another,
  title={Another look at rejection sampling through importance sampling},
  author={Chen, Yuguo},
  journal={Statistics \& probability letters},
  volume={72},
  number={4},
  pages={277--283},
  year={2005},
  publisher={Elsevier}
}

@article{gibbs2002choosing,
  title={On choosing and bounding probability metrics},
  author={Gibbs, Alison L and Su, Francis Edward},
  journal={International statistical review},
  volume={70},
  number={3},
  pages={419--435},
  year={2002},
  publisher={Wiley Online Library}
}

@incollection{bickel2008sharp,
  title={Sharp failure rates for the bootstrap particle filter in high dimensions},
  author={Bickel, Peter and Li, Bo and Bengtsson, Thomas},
  booktitle={Pushing the limits of contemporary statistics: Contributions in honor of Jayanta K. Ghosh},
  pages={318--329},
  year={2008},
  publisher={Institute of Mathematical Statistics}
}

@article{snyder2008obstacles,
  title={Obstacles to high-dimensional particle filtering},
  author={Snyder, Chris and Bengtsson, Thomas and Bickel, Peter and Anderson, Jeff},
  journal={Monthly Weather Review},
  volume={136},
  number={12},
  pages={4629--4640},
  year={2008}
}

@book{pham2009continuous,
  title={Continuous-time stochastic control and optimization with financial applications},
  author={Pham, Huy{\^e}n},
  volume={61},
  year={2009},
  publisher={Springer Science \& Business Media}
}

@article{owen2000safe,
  title={Safe and effective importance sampling},
  author={Owen, Art and Zhou, Yi},
  journal={Journal of the American Statistical Association},
  volume={95},
  number={449},
  pages={135--143},
  year={2000},
  publisher={Taylor \& Francis Group}
}

@article{hartmann2016model,
  title={Model reduction algorithms for optimal control and importance sampling of diffusions},
  author={Hartmann, Carsten and Sch{\"u}tte, Christof and Zhang, Wei},
  journal={Nonlinearity},
  volume={29},
  number={8},
  pages={2298},
  year={2016},
  publisher={IOP Publishing}
}

@article{zhang2014optimal,
	title = {Optimal control of multiscale systems using reduced-order models},
	author = {Hartmann, Carsten and Latorre, Juan C. and Pavliotis, Grigorios A. and Zhang, Wei},
	year = {2014},
	journal = {J. Computational Dynamics},
	volume = {1},
	pages = {279-306} 
}
\end{document}